\documentclass[article]{interact}

\usepackage[numbers,sort&compress]{natbib}
\bibpunct[, ]{[}{]}{,}{n}{,}{,}
\makeatletter
\def\NAT@def@citea{\def\@citea{\NAT@separator}}
\makeatother

\usepackage{color}
\usepackage{graphicx,epsfig}

\theoremstyle{plain}
\newtheorem{theorem}{Theorem}

\newtheorem{proposition}{Proposition}

\theoremstyle{definition}
\newtheorem{definition}[theorem]{Definition}

\theoremstyle{remark}
\newtheorem{remark}{Remark}

\def\vst{\hbox{\em vst\,}}
\def\sgn{\hbox{sgn}}
\def\df{\hbox{\em df\,}}
\def\ncp{\hbox{\em ncp\,}}
\newcommand{\bnu}{\bm \nu}
\newcommand\bp{\bm{p}}
\newcommand\bq{\bm{q}}
\newcommand\bu{\bm{u}}
\newcommand\bX{\mathbf{X}}
\newcommand\bx{\mathbf{x}}
\newcommand\by{\mathbf{y}}

\def\e{\hbox{E}}

\def\var{\hbox{Var}}
\def\bias{\hbox{Bias}}
\def\kld{\hbox{\em KLD\,}}\newcommand{\key}{{\mathcal{K}}}
\newcommand{\calB}{{\cal B}}
\newcommand{\calC}{{\cal C}}
\newcommand{\calS}{{\cal S}}

\begin{document}

\title{Evidence for goodness of fit in Karl Pearson chi-squared statistics}
\maketitle
\author{\name{R.~G. Staudte\thanks{Email: r.staudte@latrobe.edu.au}\\
\affil{Department of Mathematics and Statistics,  La Trobe University, Melbourne, Victoria, Australia 3086}}}

\begin{abstract}
Chi-squared tests for lack of fit are traditionally employed to find evidence
against a hypothesized model, with the model accepted if the Karl Pearson  statistic comparing
observed and expected numbers of observations falling within cells is not significantly large.
However, if one really wants evidence {\em for} goodness of fit,  it is better to adopt an
equivalence testing approach in which small values of the chi-squared
 statistic indicate evidence for the desired model. This method requires one to define
 what is meant by equivalence to the desired model, and guidelines are proposed.
  It is shown that the evidence  {for equivalence} can distinguish between normal
   and nearby models,  as well between the  Poisson and over-dispersed models. Applications
   to evaluation of random number generators and to uniformity of the digits of pi are included.
    Sample sizes required to obtain  a desired expected  evidence for
   goodness of fit are also provided.
\end{abstract}


\section{Introduction}\label{sec:introd}

 {The aim of this paper is to explain how to define and interpret evidence for equivalence to a specific
model.  We adopt the composite null and alternative hypotheses of the equivalence testing framework, and transform
the usual chi-squared distributed statistic to a normally distributed one, called the evidence for the alternative hypothesis.
We recommend interpreting this statistical evidence  to gain insight into whether one should adopt the desired model, and
provide the statistical software for  doing so. }

\subsection {Background and summary}

Tests for lack of fit based on the Karl Pearson \cite{KP-1900}  statistic have  been the subject of numerous theoretical and applied research papers, see
\cite{G-N-1996}, for example, for results and references. And, they are almost universally found in statistical textbooks,   because of their simplicity and general applicability.  The intent of the test is to validate subsequent use of  the null model, the argument being:  if the test does not reject this model at the usual levels, then it is safe to adopt the model in further analysis.

 {Critiques of this argument  appeared as early as \cite{fry-1938}, \cite{berk-1938,berk-1942}. An alternative methodology, recently proposed by Wellek \cite[Ch.8]{wellek-2003},  approaches goodness of fit by means of equivalence testing, in that what one wants to establish is placed in  the alternative  hypothesis. That is, the order of the hypotheses is reversed so that the null hypothesis that is protected is non-equivalence (the model does not fit).  This approach is not without difficulty because, to begin with, it forces one to define an equivalence model:\quad  which models are close enough to the desired one to be regarded, for all practical purposes, as \lq equivalent\rq\,?  And then  one has the  accompanying complexities of any hypothesis test:  how to choose the level; is the power of the test  for the given sample size large enough, etc.   Without in any way disparaging  the usefulness of these concepts for formal testing or analysis, this paper proposes, instead of such testing, to compute  an informative statistic  called  {\em evidence for the alternative hypothesis}.   Its practicality and foundational basis  for inference is established in  \cite{KMS-2008,KMS-2011,KMS-2014},\cite{M-S-2012}\cite{M-S-2013},\cite{PS-2016}   and in particular for equivalence testing in \cite{M-S-2016}. }

 {We proceed in steps, first finding the evidence for lack of fit in the traditional context of chi-squared tests; second, we describe Wellek's proposal for equivalence testing for goodness of fit; and third, we put both ideas together by finding the evidence for an equivalence model.
To be specific, first  in  Section~\ref{sec:evidpos}  we describe the evidence {\em for lack of fit} (and hence {\em against} the desired model) in the chi-squared statistic.   It is also shown that the expected evidence is approximately the square root of the symmetrized  Kullback-Leibler divergence  \cite{kl-1951} between the null and alternative  chi-square models. }

 {Next in  Section~\ref{sec:equivtest}  the equivalence testing approach of Wellek \cite{wellek-2003} is described, for it is the setting for the remainder of the paper. Evidence {\em for} equivalence is defined and justified in  Section~\ref{sec:equivevid}.
 It requires one to choose the boundary between equivalence and non-equivalence  and general proposals are  found  in Section~\ref{sec:choosingbdy}. This leads to a formula for choosing the minimum sample size required to obtain a desired degree of expected evidence in Section~\ref{sec:choosingn}.
Examples follow in Section~\ref{sec:exsevid}, including  evidence for normality over nearby models and  evidence for the Poisson model in the presence of over-dispersion.
Evidence for uniformity of digits produced by a random number generator and the decimal digits of $\pi$ are relegated to the on-line supplementary material.  Numerous related research topics are proposed in Section~\ref{sec:summary}.}

\subsection{The traditional Karl Pearson lack of fit test}\label{sec:defnsKP}

 Many lack of fit tests are based on test statistics $S$
having an approximate central or non-central chi-squared distribution with known degrees of freedom (\df) $\nu $ and unknown non-centrality parameter (\ncp) $\lambda \geq 0,$ see e.g. \cite{G-N-1996}; this assumption is abbreviated  $S\sim \chi ^2_{\nu ,\lambda} .$  The Karl Pearson
lack of fit test arises as follows: given a sequence of independent trials indexed by $k=1,2,\dots ,n$ with outcomes lying  in one of $r$ mutually exclusive sets (cells) and with respective probabilities $p_1,\dots ,p_r$, let $\nu _i $ be the frequency of outcomes of cell $i$ in the first $n$ trials, $i=1,\dots ,r$. Then $\bnu =(\nu_1,\dots ,\nu_r)^\top$ is a sufficient statistic for $\bp =(p_1,\dots ,p_r)^\top$ and $\bnu $ has a multinomial distribution ${\cal M}(n,\bp )$, see for example  \cite[pp.1--3]{G-N-1996}.
Let \[\bX _n = \left( \frac{\nu _1-np_1}{\sqrt{np_1}},\dots, \frac{\nu _r-np_r}{\sqrt{np_r}}   \right)^\top~.\]
The \citep{KP-1900} chi-squared statistic is based on the length squared of this vector, namely
\begin{equation}\label{eqn:KP}
   \| \bX _n\| ^2= \bX _n^\top\bX _n = \sum _{i=1}^r\frac{(\nu _i-np_i)^2}{np_i}~.
\end{equation}
This  statistic $\| \bX _n\| ^2$ is sometimes simply written $\bX ^2_n$.  The hypothesized model $\bp$ is rejected if the
statistic $S=\bX ^2_n$ exceeds a pre-chosen critical point $c_{n,\alpha }$.  For large $n$, and a fixed level $\alpha $, it is known that $\Pr \{S\geq c_{n,\alpha }\}\approx \alpha $ provided $c_{n,\alpha }=\chi ^2_{r-1}(1-\alpha )$ the $1-\alpha $ quantile of the central chi-squared distribution with $r-1$ \df.  This defines the chi-squared test of Karl Pearson \cite{KP-1900}.  For the subsequent historical development of this famous test, see the references in \cite[p.7]{G-N-1996}.  Further, under certain alternatives (\ref{eqn:lambda}) to the null the distribution of $S$ is well-approximated by the non-central chi-squared distribution   $ \chi ^2_{\nu ,\lambda} $, where $\nu =r-1$ and the non-centrality parameter (\ncp ) \;$\lambda >0$.
  {Thus the Karl Pearson test for lack of fit based on $S$, is essentially testing the simple null hypothesis $\lambda =0$ against the composite alternative $\lambda >0.$} The asymptotic power function of such  a level-$\alpha $ test for lack of fit is then given by
\begin{equation}\label{eqn:asympowpos}
  \Pi  _\alpha (\lambda )=P\{S\geq  c\}, \hbox { where }  c=\chi ^2_{\nu,\lambda }(1-\alpha )~.
\end{equation}

\subsection{Evidence for lack of fit in the  Karl Pearson statistic}\label{sec:evidpos}

Long ago the medical researcher Berkson \cite{berk-1942} took issue with null hypothesis significance testing. He claimed:  { \em \lq Nor do you find experimentalists typically engaged in disproving things. They are looking for appropriate evidence for affirmative conclusions.\rq }\   Such  appropriate evidence  is
  often routinely found by statisticians while carrying out a test.  To be specific, consider the evidence for an alternative hypothesis as introduced in \cite{KMS-2008}.   In its simplest context, one has data $X$ normally distributed with mean $\theta $, variance 1, denoted  $X\sim N(\theta ,1)$,  and wants information regarding a null  hypothesis $\theta \leq \theta _0 $  and alternative $\theta >\theta _0$. The {\em evidence for the alternative} is then defined to be  $T=X-\theta _0$.  This evidence is an estimator of its mean with a  standard normal error, and so is written $T\pm 1$.
 Values of $T$ near 1.645,  3.3 and 5 were suggested by \cite[p.17]{KMS-2008} as \lq weak\rq,\  \lq moderate\rq\  and \lq strong\rq\  evidence for the alternative, because for a level $\alpha $ test of the above hypotheses having power $1-\beta(\theta)$  the expected evidence satisfies
\begin{equation}\label{evidlevelpower}
  \e _\theta [T]= \Phi ^{-1}(1-\alpha )+\Phi ^{-1}(1-\beta(\theta) )~.
\end{equation}
Here $\Phi ^{-1}$ is the inverse of $\Phi $,  the standard normal distribution function.
Thus a level-0.05 test based on $T$ has power  1/2, 0.95, or 0.9996  when the expected evidence for the alternative is respectively  weak (1.645),  moderate (3.3) or strong (5.0) .

There are other reasons for adopting these rough descriptive labels. Usually $T=T_n$ has expected value growing with the sample size $n$ at the rate $\sqrt n\,.$  Moderate evidence for an alternative  can be expected if an experiment is repeated under the same conditions as one that yielded weak evidence for it, provided the sample size is quadrupled. In symbols, having found $T_{n}=1.645 \pm 1$  one can expect in a replicated experiment the test statistic $T^\prime_{4n}=3.3\pm 1$ .
Also, when choosing sample sizes one often stipulates power 0.8 at level 0.05 for a specific  alternative, which is considered a minimal requirement. Such a
test will have a   expected evidence of 2.5,  which is about  halfway between weak and moderate. Studying evidence for the alternative hypothesis  is in this respect more basic than examining its contributing elements level and power, one of which is open to arbitrary choice.

Negative evidence for the alternative can be interpreted as positive evidence for the null hypothesis.  For example, $T=-3.3\pm1$ is moderate evidence for the null.  This feature enables \lq non-significant\rq\  results  to be  easily combined with \lq significant\rq\  ones in a meta-analysis, see \cite{KMS-2008,KMS-2011,KMS-2014} for examples.   All the above  advantages of evidence for the alternative hypothesis justified
its introduction  as a useful  {\em operational} one. It works!

  {Why it works was explained when a
 {\em foundational} rationale for this  calibration scale was described in  \cite{M-S-2012,M-S-2013,M-S-2016}, for exponential families, among others;
  and for the difference of two proportions in \cite{prst-2014}.
   The examples in these papers suggest that quite generally the evidence for the alternative has an asymptotic mean that  is approximately the signed square root of the symmetrized Kullback-Leibler divergence between null  and alternative distributions and further it has a standard normal error. The evidence for the alternative
    is  an effect size on a simple canonical scale, where the unit of measurement is the standard error of the effect size.  One usually transforms a test statistic $S$ to this simple normal calibration scale with  a variance stabilizing transformation (\vst), say $T=T(S)$.  A specific example is given below in (\ref{eqn:evidpos}),(\ref{eqn:biasadjevidpos})
   and plotted in Figure~\ref{fig2}(a).}

\begin{remark}\label{rem1}
 Traditionally, given an observed $S=s$ and its  {transformed} value $T(s)=t$  one would compute a p-value $PV(t)=1-\Phi(t)$, but this has the disadvantage of moving from a useful calibration scale for evidence for the alternative  to the  p-value scale where the result is well known to be open to  misinterpretation.  For a thorough discussion, see \cite[pp. 4-5, 113-119]{KMS-2008} or \cite{wass-2019}. Another advantage of this calibration scale is that it  reminds the user that statistical outcomes are subject to error, in this case a unit normal error, while the p-value is reported to two or more decimal places,  giving the impression of precision while hiding information about the randomness that yielded the result.
\end{remark}

Let $S\sim \chi^2_{\nu ,\lambda }$  be a statistic for testing $\lambda =0$ against $\lambda >0.$  The mean and variance of
$S$ are $\e_{\nu ,\lambda }[S]=\nu +\lambda $ and  $\var _{\nu ,\lambda }[S]=2\nu +4\lambda .$
The basic idea is to transform the test statistic $S$ into $T=T(X)\sim N(\e_\lambda [T],1)$  with $\e_\lambda [T]$ increasing from 0 as the parameter $\lambda $ moves away from the null. A standard derivation of the \vst\,  based on the delta method, see \cite[p. 183]{KMS-2008}, for example,
yields   $T_1(S)= \sqrt {S-\nu /2}+c_1 $, where  $c_1$ is an arbitrary constant.     {The authors of \cite{KMS-2008}  note that  this derivation
is only valid for $S \geq \nu $,  and one needs to smoothly extend the domain of $T$  to $0\leq S< \nu $ so that evidence $T(S)$ is defined for all $S\geq 0$ and strictly increasing (so that it is a test statistic). They
  suggested a \vst\, extension for all $S \geq 0$ by means of a symmetrization argument about  the median.  Here is proposed a less complicated solution, as follows.}

\begin{figure}[t!]
\centering
\includegraphics[width =10cm,height=8cm]{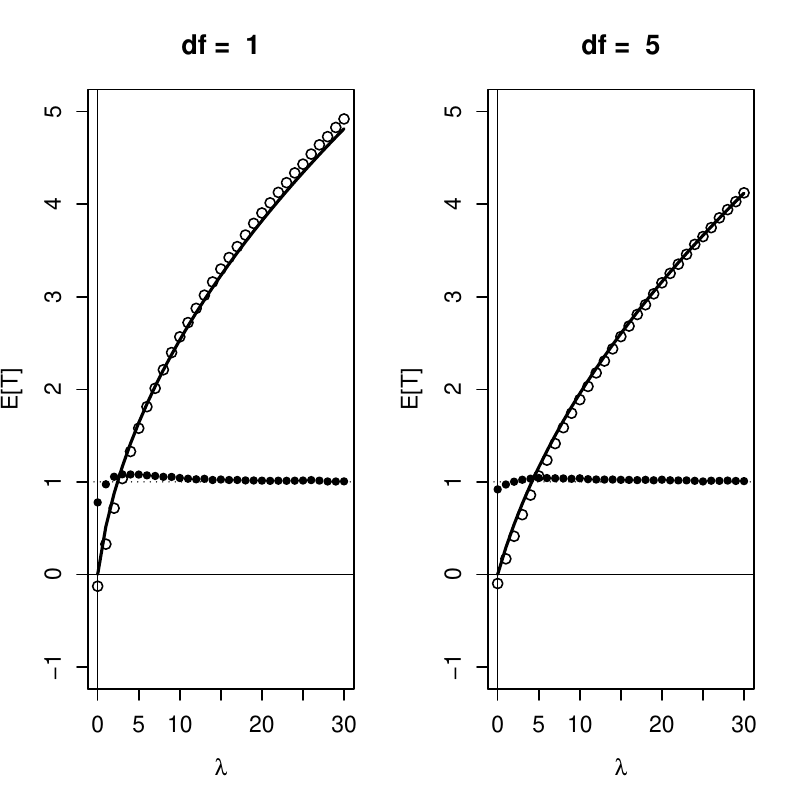}
\caption{\small The left-hand  plot shows the asymptotic mean (\ref{eqn:ETpos})  plotted as a function of $\lambda $ for $\df=\nu =1$. The circles and dots are respective simulated values of $\bar T$ and $ s_T$,  the mean and standard deviation of 40,000 replications of $T(S)$, for $S \sim  \chi  ^2_{1,\lambda}$,  for $\lambda  = 0:30/1. $ }       \label{fig1}
\end{figure}

\begin{figure}
\centering
\includegraphics[width =10cm,height=8cm]{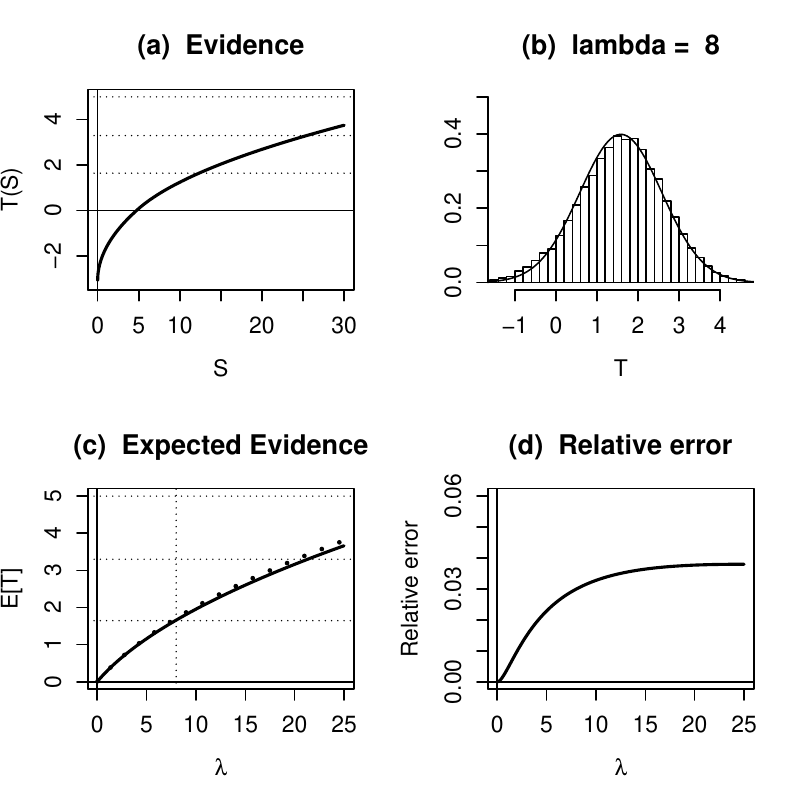}
\caption{\small  For $\nu =5$ throughout, plot (a) shows the evidence (\ref{eqn:biasadjevidpos}) as a function of $S$.
  Plot (b) is a histogram of simulated values of   $T_{ba}(S)$, for  $S\sim \chi ^2_{5,8},$ with a superimposed normal density. In Plot (c) the asymptotic mean (\ref{eqn:ETpos})  is plotted as a function of $\lambda $ as a solid line, while the dotted  line shows $\sqrt {J(\lambda ,\lambda _0)}\,$ with the relative error of this approximation given in Plot (d).  }\label{fig2}
\end{figure}

For  the central chi-squared distribution $\lambda =0$ the transformation $T_0 (S)=\sqrt {2S }\ +c_0$ has variance near one for all $\nu >0$ and one can expect   this also to be true in the non-central case for small  $\lambda $. Further, by  choosing $c_0= -\sqrt{2\nu }$ the expected value of $T_0(S)$ should be near 0 for small $\lambda$.  By piecing together the transformations $T_0,T_1$ one obtains a  \vst \, $T=T(S)$ that has a positive derivative for all $S>0$ and satisfies $T(\nu )=0.$
 \begin{equation} \label{eqn:evidpos}
T=T(S) = \left\{
             \begin{array}{ll}
               T_0(S)=\sqrt {2S }-\sqrt{2\nu }\, & 0\leq S <  \nu  ; \\
               T_1(S)=\sqrt {S-\nu /2}-\sqrt{\nu /2} , & \nu \leq S .
             \end{array}
           \right.
\end{equation}
For fixed $\nu \geq 1$  and  $\lambda >0$ one expects that  this $T=T(S)\sim N(\e _{\nu,\lambda } [T],1).$
As explained in  Section~\ref{app:biaspos}, to first order $ \e _{\nu,\lambda } [T]$ is simply given by
\begin{equation}\label{eqn:ETpos}
  \e _{\nu,\lambda }[T] \doteq   \sqrt{\lambda+\nu/2}\,-\sqrt{\nu/2} \,~.
\end{equation}
It is also shown there that $T$ has a  negative bias for (\ref{eqn:ETpos}), which leads to a bias-adjusted version:
\begin{equation}\label{eqn:biasadjevidpos}
  T_\text {ba}(S) = T(S)+0.2/\sqrt \nu\,~, \hbox { for all $S\geq 0$~ .}
\end{equation}
 Simulations confirm that, to a good approximation,  $T_\text {ba}\sim N(\e _{\nu,\lambda } [T],1)$;  these simulations   can be carried out for various $\nu $ and $\lambda $ using R software \citep{R} scripts in the on-line supplementary materials.  Some (typical) results are shown in Figures~\ref{fig1} and \ref{fig2}.
The plots in Figure~\ref{fig1} demonstrate that the approximation (\ref{eqn:ETpos}) is quite good for $\nu =1$ and $\nu =5$ for
values of $\lambda \geq 0 $ of interest.  Further the standard deviations of the simulated $T_\text {ba}(S)$ values are quite near one.

 {The results of another  study for $\nu =5$ are depicted in  Figure~\ref{fig2}. The graph of the evidence for $\lambda >0$
is shown in plot (a), and to be compared with the horizontal dotted lines, which are rough guides to what is considered weak, moderate and strong evidence for $\lambda >0$. }

 {In Plot~(b) of Figure~\ref{fig2} a histogram summarizes 40,000 simulated values from $T_\text {ba}(S), $ with $S\sim \chi ^2_{5,8},$  where $\lambda =8$ is a possible alternative.  Note that the histogram is close to a normal density with standard deviation one, as expected.   Further experimentation with $\nu =1$ or 2 and $\lambda $  near 0 reveals  that in these cases a slightly truncated normal distribution results but  this shortcoming does not materially affect applications.}

 {Plot~(c) of Figure~\ref{fig2} shows  the expected evidence (\ref{eqn:ETpos}) against $\lambda =0$ as a solid line; it has value 1.65  at $\lambda =8$, marked by a dotted vertical line.  Thus on average, when $\lambda =8$ there is weak evidence for the alternative $\lambda >0.$
Also in Plot~(c) is shown the graph (dotted line) of
the numerically computed  square root of $J$, where $J$ is the Kullback-Leibler  symmetrized divergence  \cite{kl-1951}  between the densities of $\chi ^2_5$ and $\chi ^2_{5,\lambda };$ it is very close to the asymptotic expected evidence (\ref{eqn:ETpos}) for a wide range of $\lambda .$}

\subsection{Example 1:  Evidence for biasedness of a die}\label{sec:wellekex}
 {The difference between the classical testing for lack of fit described in Section~\ref{sec:defnsKP} and the evidence for the alternative of lack of fit
 are now illustrated on the die tossing  example   of \cite[Ex.8.1]{wellek-2003}.}
   A die with sides numbered $(1,2,3,4,5,6)$ is tossed $n=100$ times with resulting counts $\bnu = (17,16,25,9,16,17)$.
In the traditional test for biasedness  the null hypothesis of unbiasedness  is rejected if the  Karl Pearson statistic (\ref{eqn:KP}) is \lq too large\rq.\  For these data $S=7.76$, which leads to a p-value computed from the central chi-squared distribution $P(\chi ^2_5\geq 7.76)\approx 0.17.$ Because this  is not significant
at the usual levels, a decision would usually be  made to accept unbiasedness. Wellek \cite{wellek-2003} queries such a decision and  then uses equivalence testing to show that this decision is unwarranted, see Section~\ref{sec:eqtestwellek}.

 {It is of interest to find the evidence for biasedness in these same data using the normal calibration scale.  The  Karl Pearson statistic for $n=100$ is approximately distributed as  $\chi ^2_{5,\lambda }$ for some  unknown $\lambda \geq 0$, which is 0 if the die is unbiased, and otherwise positive. Here $\nu =r-1=5$ and from (\ref{eqn:biasadjevidpos})  the evidence for biasedness $\lambda >0$ is $T_\text {ba}(S)=\sqrt {S-\nu /2}-\sqrt{\nu/2}+0.2/\sqrt{\nu }\,=\sqrt{7.76-2.5}\, -\sqrt{2.5}\,+0.2/\sqrt{5}=0.8, $ about halfway between 0 (negligible) and 1.645 (weak) evidence. This estimate of the unknown expected evidence has standard error 1 and so there is positive, but very weak, evidence for biasedness as measured by $T_\text {ba}(S)=0.8\pm1 $  in the results  $\bnu = (17,16,25,9,16,17)$ of the $n=100$ tosses of the die.  {\em  Clearly, a conclusion of biasedness is not warranted by these data.}
 If the outcomes had been   $\bnu = (23,10,25,9,16,17),$ for example, then the reader can check that
$S=12.8$ and $T  _\text {ba}(S)=1.72\pm 1 $, which is weak evidence for biasedness.}

 {It is possible to obtain negative evidence for the alternative of biasedness in the above 100 tosses of a die.  For an extreme example,  if $\bnu $ were $(17,17,17,16,16,17)$, then $S=0.08$ and $T  _\text {ba}(S)=-2.67$.  Because this  is negative (weak to moderate) evidence for biasedness $\lambda >0$,  it can be interpreted as  positive evidence  $2.67\pm 1 $   for unbiasedness $\lambda =0$.  However, this observation is of limited practical use, because one is very unlikely to obtain a value of $S$ so near 0 in 100 tosses, even if the die were perfectly fair.  Further reflection reveals that for these hypotheses and sample size, one is unlikely to obtain even weak evidence by this method for unbiasedness when the die is perfectly fair.  One needs a larger sample size or larger hypothesis containing the desired model. In this connection, Greenwood and Nikulin \cite[p.27-28]{G-N-1996}, for example,  point out that a simple null hypothesis is never exactly true and so it would be better to test the composite null $0\leq \lambda \leq \lambda _0$ against the composite alternative $\lambda >\lambda _0$ for some given positive $\lambda _0$.  For the uniform case just discussed they recommend $\lambda _0=1/n.$  This approach is rarely carried out in practice, perhaps because it is complicated:  one then needs the {\em non-central} chi-squared distribution to find the critical point; or, more likely, one is still uncertain about how to choose $\lambda _0$.  Further, interpreting negative evidence for an alternative as positive evidence for a null hypothesis is an awkward way to formulate the problem we are trying to solve.
Rather than continue down this path,  we now move to the equivalence testing framework, which is more natural for finding evidence for the desired model or an equivalent one. }

\section{The equivalence testing approach to goodness of fit}\label{sec:equivtest}
The  equivalence testing approach   {of  Wellek \cite{wellek-2003}
 to establishing goodness of fit is a natural approach to the problem.  However, his normal approximation theory is somewhat complicated,
 so in Section~\ref{sec:eqtestchisq}  we describe a parallel equivalence testing method for goodness of fit based on the chi-squared approximation. }

\subsection{Wellek's proposal for establishing goodness of fit}\label{sec:eqtestwellek}
Recall from Section~\ref{sec:defnsKP} that the frequency vector  $\bnu $ for the $r$ cells has a multinomial distribution ${\cal M}(n,\bp )$.
 {Wellek \cite[Sec.8.1]{wellek-2003}  proposes that for comparing two multinomial distributions ${\cal M}(n,\bp )$  and  ${\cal M}(n,\bp ^0)$
one use the  squared Euclidean distance $d^2=d^2(\bp,\bp _0)=\sum _i(p_i-p_i^0)^2$ . } The null hypothesis of non-equivalence is stated as $d^2\geq d_0^2$ for some fixed boundary $d_0$, while the equivalence alternative is $0\leq d^2< d_0^2$; in both cases $\bp $ varies over the $(r-1)$-simplex in $r$ dimensional space.  Letting $\hat \bp =\bnu /n$ denote the
maximum likelihood estimator of $\bp $, \cite{wellek-2003} proceeds to derive the asymptotic  distribution of
$\hat d^2= d^2(\hat \bp, \bp ^0) $, including an expression for the asymptotic variance in terms of $\bp ^0 $ and $\bp$ which must
be estimated.    His equivalence test  for goodness of fit rejects the null  $d^2\geq d_0^2$ (lack of fit)
at level $\alpha $ in favour of equivalence if $\hat d^2(\hat \bp, \bp ^0) $ is smaller than the $\alpha $ quantile of the approximating
normal distribution. For the example of Section~\ref{sec:wellekex}, where the ideal model is uniform $  \bp ^0=\bu=(1/6,\dots ,1/6)$ and $n=100$,  \cite{wellek-2003}  chooses $d_0^2=0.15^2$.  The data  $\bnu = (17,16,25,9,16,17)$ gives $\hat \bp =\bnu/n$, which leads to
 $\hat d^2= d^2(\hat \bp, \bp ^0)=0.02913 $ which is not significant at level $\alpha =0.05$. He concludes:
\begin{quote}\em Thus the example gives a concrete illustration of the basic general fact (obvious enough from a theoretical viewpoint) that the
traditional $\chi ^2$ goodness of fit to a fully specified multinomial distribution is inappropriate  for  {\em establishing} the hypothesis
of (approximate) fit of the true to the pre-specified distribution.
\end{quote}

\begin{figure}[t!]
\centering
\includegraphics[scale=0.8]{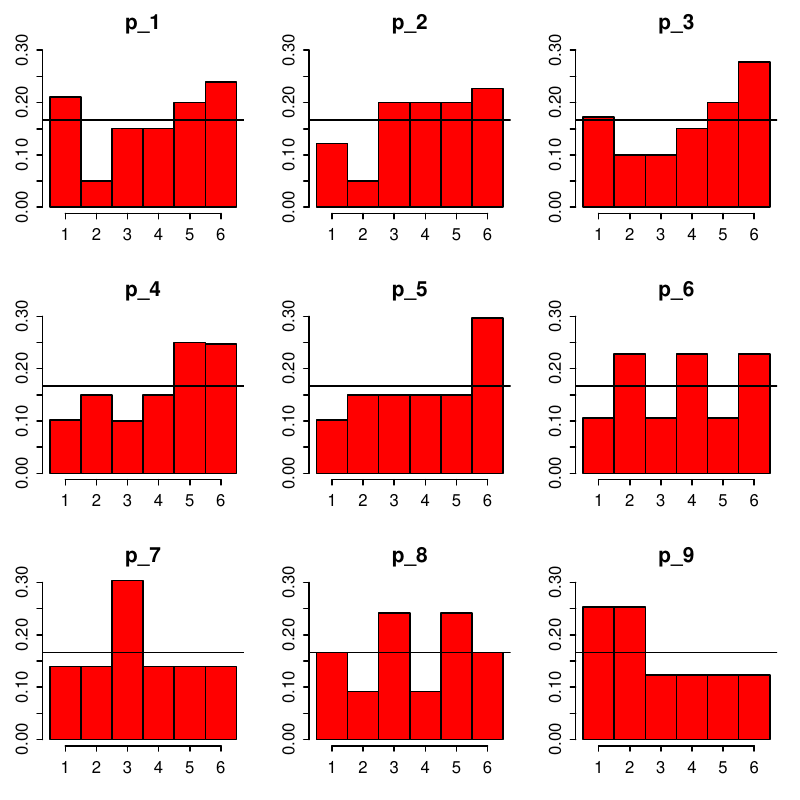}
\caption{All 9 models have Euclidean distance $d_0=0.15$ from the uniform.   The top six are taken from \cite[Table 8.2b]{wellek-2003}.
His choices are supplemented with  another three models $\bp _7,\bp_8,\bp _9$  all of distance 0.15 from the uniform, for the reader's consideration.
Which of the nine models is closest to uniformity? Which is furthest?}  \label{fig3}
\end{figure}

\noindent Wellek makes  another  point about the above example:  the {\em level} $\alpha =0.05 $ based on the asymptotic normal distribution
when $d^2=d_0^2=0.15^2$ does not necessarily describe the {\em size} of the test.  He lists six models $\bp _1,\dots ,\bp_6$
for which   $d(\bp_i,\bu) =0.15$ to five decimal places and are plotted here in the top two rows of Figure~\ref{fig3}. His Table 8.2b gives the actual sizes of the tests;  they range from 0.00833 for $\bp _2$ to 0.03943 for $\bp _5$.

These plots raise the question of what is meant by   {\lq equivalence to uniformity\rq.\ } One can measure non-uniformity of  $\bp =(p_1,\dots, p_r)$ from  $\bu _r=(1/r,\dots ,1/r)$  by the Euclidean distance $d(\bp ,\bu_r)$, the sup metric $M(\bp ,\bu _r)=\max _i\{|p_i-1/r| \}$ or by a semi-metric such as
 the symmetrized Kullback-Leibler divergence,   {shown in Appendix~\ref{app:kld} to be}  $J(\bp ,\bu_r)=\sum _i (p_i-1/r)\ln(p_i)$.

 Table~\ref{table1} lists these values for
the models of Figure~\ref{fig3}.    Model $\bp _6$ is closest in the sup metric while model $\bp _7$ is furthest.  Model  $\bp _2$ has the largest divergence, while $\bp _7$ the smallest.
None of the models in Figure~\ref{fig3} would likely be considered \lq equivalent\rq\  to the uniform model in practice; a bound on $M$ such as $M\leq M_0= k/r$  seems desirable.  This allows  for a $100k$\% relative error in each probability.
 It is clear that  model $\bp _7$ in Figure~\ref{fig3} has the highest possible relative error $k=0.82$ subject to $d_0=0.15$.
  For  recommendations for choosing $d_0$, see Proposition~\ref{prop1} and accompanying remarks in
Section~\ref{sec:choosingbdy}.

\begin{table}[b]
 \tbl{Nine possible models which are distance $d=0.15$ from uniformity are plotted in Figure~\ref{fig3}.  The sup metric distance $M$
and symmetrized Kullback-Leibler divergence $J$ of $\bp$ from $\bu $ is tabled for each model.  Also listed is the estimated power
$\Pi  _{0.05} (\bp _j)=P_{\bp_j}\{S\geq  c\}$, where $c=\chi ^2_5(0.95)$ of the level 0.05 $\chi ^2$ test for non-uniformity. }
  {\begin{tabular}{cccccccccc} \toprule
       Model                        &  $\bp_1$   &$\bp_2$ & $\bp_3$ &$\bp_4$& $\bp_5$ &$\bp_6$ & $\bp_7$ &$\bp_8$ &$\bp_9 $ \\
       \midrule
         $d$                           &  0.150    &  0.150  &  0.150  & 0.150   & 0.150   &  0.150    &   0.150    & 0.150   &  0.150    \\
        $M$                          &  0.117     &  0.117  &  0.110  & 0.083   & 0.131   &   0.061   &   0.137   &  0.075   &   0.087  \\
         $J$                           &  0.187     &   0.191 &  0.132  & 0.134   & 0.115   &   0.142   &   0.107   & 0.145   &   0.125  \\
 $\Pi _{0.05}(\bp _j)  $    &  0.887    &  0.899  &  0.820  & 0.823   & 0.779   &   0.836   &   0.762   & 0.839   &   0.798  \\
\bottomrule
  \end{tabular}} \label{table1}
\end{table}

\subsection{Equivalence testing with the chi-squared statistic }\label{sec:eqtestchisq}
       Hereafter it is  assumed $n$ is large  enough so that the Karl Pearson statistic for comparing these multinomial distributions has an approximate
 $\chi ^2_{\nu ,\lambda }$ distribution. An  {\em equivalence boundary value}  $\lambda _0$ that separates the null hypothesis of non-equivalence
 $\lambda \geq \lambda _0$ from the equivalence alternative $0\leq  \lambda <\lambda _0  $  must be chosen in advance.
Then one can  carry out a test rejecting non-equivalence  at level $\alpha $ when $S\leq c_\alpha=\chi^2_{\nu \lambda _0, } (\alpha )$, the $\alpha $-quantile of the $\chi ^2_{\nu ,\lambda _0}$ distribution.

The asymptotic distribution of Pearson's statistic under alternatives plays an essential role in what follows. The alternative hypotheses are denoted for each $n$ by $\bp ^{(n)}=(p_1^{(n)},\dots ,p_r^{(n)} )$. The basic additional assumption is that for some $\lambda >0$
\begin{equation}\label{eqn:lambda}
   \lambda _n=  n\sum _{i=1}^r\frac{(p _i-p_i^{(n)})^2}{p_i}\to \lambda ~ .
\end{equation}
Then as $n\to \infty $, under the sequence $\{\bp ^{(n)}\}$ Pearson's chi-squared statistic (\ref{eqn:KP}) has a non-central chi-squared distribution in the limit with \df\; $\nu =r-1$ and \ncp \ equal to $\lambda $; see  \cite[Sec. 3]{G-N-1996}.  Using this result, one can find an approximate power function of an asymptotic level-$\alpha $ chi-squared test.  Letting $c_\alpha=\chi^2_{\nu ,\lambda _0}(\alpha )$, it is
$\Pi _\alpha (\lambda ) = P( \chi^2_{\nu ,\lambda }     \leq c_\alpha )$  for $ 0\leq \lambda <\lambda _0.$
 In particular for  $\bp =\bu _r =(1/r,1/r,\dots ,1/r) $  one obtains $\lambda _n =n\,r\,d_n^2,$  where $d_n=d(\bp ^{(n)} ,\bu _r)$.
Thus in a Euclidean neighborhood of the uniform, the Karl Pearson statistic has an approximate $\chi ^2_{\nu ,\lambda _n}$  distribution,
with $\nu =r-1$ and $\lambda _n =rn\, d_n^2$.

\begin{remark}\label{rem2}  { In applications, one does not always have  sequences $\{\bp ^{(n)}\}$, $\{\lambda _n\}$ with a limit $\lambda $ in mind, but assumes that  the observed
$\hat {\bp }^{(n)}$ and corresponding $\hat {\lambda }_n$ could so be embedded, with  $\lambda =\hat {\lambda }_n$.
Such license is  subject to  rules of thumb such as $n/r \geq 5;$ and then it is assumed that the  $\chi ^2_{\nu ,\lambda _n}$
distribution is a good one to approximate the  distribution of the Karl Pearson statistic, provided $\nu =r-1-s$ and $s$ is the number of estimated model parameters.
Further many of these limit theorem proofs allow the number of cells $r$ to be increasing with $n$, but only at a much smaller rate such as $r(n)=\ln (n)$.
Caution is necessary for small and even moderate $n$ and has been widely discussed in the literature, see the text and references in \cite{cochran-1952}, \cite[pp. 18-21]{G-N-1996} and \cite{KV-2018}.  We have not encountered any problems with the accuracy of chi-squared approximations to the Karl Pearson statistic for the choices of $n$ and $r$ in our study of
evidence for equivalence. Generally speaking, equivalence testing requires larger sample sizes than needed for routine one-sided testing
and the same is true when finding evidence for  the alternative in this setting; see \cite{M-S-2016} for examples. If the chi-squared approximation to the Karl Pearson statistic
is justified in the usual setting it will also be good enough for the larger sample sizes of equivalence testing.}
\end{remark}

\section{Evidence for  equivalence  in  chi-squared statistics}\label{sec:equivevid}

\subsection{Defining the evidence for equivalence}\label{sec:defnevid}

 {As in Section~\ref{sec:eqtestchisq},  $\lambda _0$  denotes the  \lq equivalence boundary value\rq\ that separates the null hypothesis of non-equivalence
 $\lambda \geq \lambda _0$ from the equivalence alternative $0\leq  \lambda <\lambda _0  $.
Proceeding as in Section~\ref{sec:evidpos}  with $S$ the statistic (\ref{eqn:KP}) but now taking into account the reversal of hypotheses,
 transform $S$ to the normal calibration scale.  A  transformation to evidence   is a continuous,{\em strictly decreasing}
function $h(S)$ of the test statistic $S$, which for all values of $\lambda$ has an approximate normal distribution with variance one
and whose  mean function is increasing from 0 at $\lambda =\lambda _0$ as $\lambda $ {\em decreases}.}

As in Section~\ref{sec:evidpos} compose two \vst s, one for each of the regions $0\leq S<\nu$ and $\nu \leq S$:
Let $c_1=\sqrt{\lambda_0+\nu/2}\;$, and  $c_0=c_1-\sqrt{\nu/2}\, +\sqrt {2\nu }\,$  in the following expression:
\begin{equation}\label{eqn:evidchisqeq}
    T_{\lambda _0}(S) = \left\{
            \begin{array}{ll}
          c_0\, -\sqrt {2S}   & \hbox{\quad for $0\leq S < \nu $\,;} \\
          c_1\,- \sqrt {S -\nu/2}  & \hbox{\quad for $\nu \leq S $\,.}
            \end{array}
          \right.
  \end{equation}
This $T_{\lambda _0}$ is not only a differentiable, strictly decreasing function, it is a transformation to evidence  for equivalence that, to first order, has expectation:
  \begin{equation}\label{eqn:keychisq}
  \e _{\nu ,\lambda }[T_{\lambda _0}(S)]\doteq\,  \key _{\nu,\lambda _0}(\lambda ) \equiv \sqrt {\lambda _0+\nu /2}-\sqrt {\lambda +\nu /2}~ \hbox { for }
0\leq \lambda \leq \lambda _0~.
\end{equation}

\begin{figure}
\centering
\includegraphics[width =10cm,height=8cm]{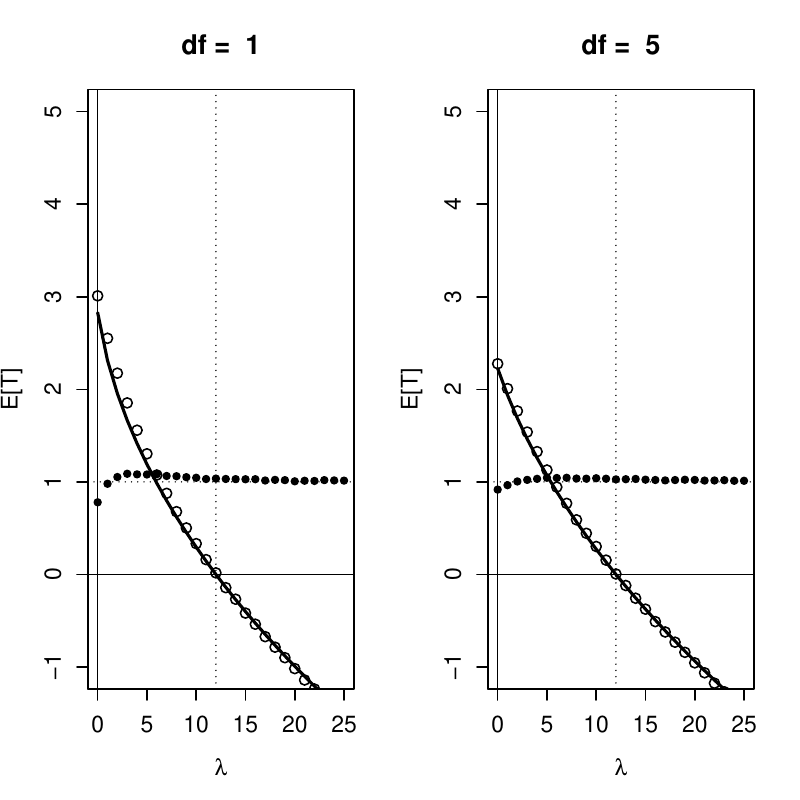}
\caption{\small The left-hand  plot shows the asymptotic mean (\ref{eqn:keychisq}) for $\df=\nu =1$.  The circles are simulated values of $\bar T_\text{ba}$,  the mean of 40,000 replications from $T_{\lambda _0,\text {ba}}(S)$, for $S \sim  \chi  ^2_{1,\lambda }$ ,  for $\lambda  = 0:25/1. $    The black dots are the standard deviations of these samples. }  \label{fig4}
\end{figure}

The evidence for equivalence $T_{\lambda _0}(S)$ defined by (\ref{eqn:evidchisqeq}) has a non-trivial upward bias $ 1/(2\sqrt{\lambda+\nu/2}\,).$, as shown in Appendix~\ref{app:bias0}.
 While $\lambda $ is unknown,
one can remove the bias at the boundary point $\lambda _0$ (and, it turns out,  smaller $\lambda $) by defining:
\begin{equation}\label{eqn:biasadjevideq}
  T_{\lambda _0,\text {ba}}(S) = T_{\lambda _0}(S)-1/(2\sqrt{\lambda_0+\nu/2}\,)~, \hbox { for all $S$ .}
\end{equation}
 This bias adjusted  $T_{\lambda _0,\text {ba}}(S) $  has expectation very near  the asymptotic mean $\key _{\nu,\lambda _0}(\lambda )\, $  defined in (\ref{eqn:keychisq}) for $\nu \geq 1$ and the region of interest $0\leq \lambda  \leq \lambda _0$.
Examples are shown in Figure~\ref{fig4} for $\lambda _0=12.$   Not only are the biases small, but the standard deviations of the transformed
values are near one.

\begin{figure}
\centering
\includegraphics[width=10cm,height=7cm]{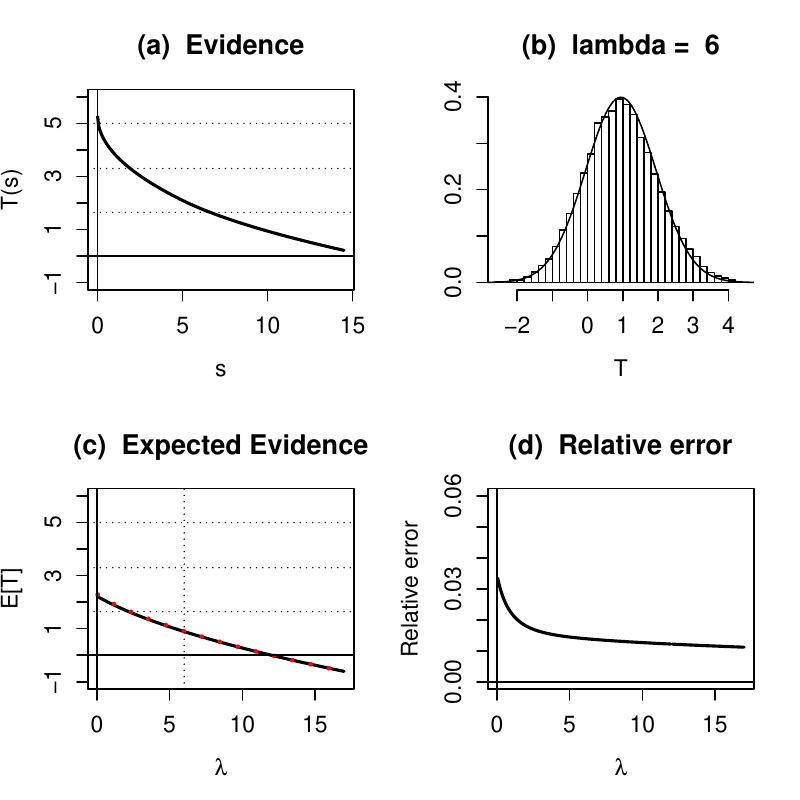}
\caption{\small  Plot (a) is the  graph of evidence $T_{\lambda _0,\text {ba}}(S)$ defined in (\ref{eqn:biasadjevideq}) when $\nu =5$ and $\lambda _0=12.$  Plot(b) is a histogram of 40,000 $ T_{\lambda _0,\text {ba}}(S)$ values,  where $S\sim\chi ^2_{5 ,6}$. Plot (c) shows the asymptotic mean (\ref{eqn:keychisq})  of $T_{\lambda _0,\text {ba}}$  plotted as a  as a solid line, while the dotted  line is the  $\sgn (\lambda _0-\lambda ) \sqrt {J(\lambda _0,\lambda )}\,$.} \label{fig5}
\end{figure}

An example of the transformation (\ref{eqn:biasadjevideq}) to evidence for the equivalence hypothesis  when $\nu =5$ and $\lambda _0 =12$  is shown in Plot(a) of Figure~\ref{fig5}.  Note that $T_{\lambda _0,\text {ba}}(S) $ is a smooth decreasing function of $S$ with maximum possible value
$T_{\lambda _0,\text {ba}}(0) =\sqrt{\lambda_0+\nu/2}\;-\sqrt{\nu/2}\, +\sqrt {2\nu }\,-1/(2\sqrt{\lambda_0+\nu/2}\,)=5.26.\,$ Thus with these parameters it is possible to get strong evidence for equivalence.  In the next plot observe
 what happens if indeed $\lambda =\lambda _1=6$.

 In Plot (b) of Figure~\ref{fig5} is shown  the histogram of 40,000 random $T_{\lambda _0,\text {ba}}$ values, obtained from random chi-squared values when $\nu =5,\lambda _0=12$ at a specific $\lambda_1 =6$. For these parameter values, the asymptotic mean (\ref{eqn:keychisq}) for equivalence is $\key _{5,12}(6)=\sqrt {12+2.5}-\sqrt {6 +2.5}\,=0.89,$ which is very weak. The sample mean and standard deviation of these
$T_\text {ba}$ values are respectively $\bar T_{\lambda _0,\text {ba}}=0.95$ and $s_T=1.03.$ In general, the asymptotic mean (\ref{eqn:keychisq})
of $T_{\lambda _0,\text {ba}}$ and its expected value $\e _{\nu, \lambda }[T_{\lambda _0,\text {ba}}]$ are very close for parameters of interest.
 The superimposed normal density with these parameters suggests that it is also approximately normal.

 {An anonymous referee pointed out the strong resemblance between the formulae in (\ref{eqn:evidchisqeq}) and (\ref{eqn:evidpos}); and it is true that if one substitutes $\lambda _0=0$ into the formulae (\ref{eqn:evidchisqeq}), then $T_{\lambda _0}(S)$ is indeed the negative of $T(S)$ defined by (\ref{eqn:evidpos}). However,
 $\lambda _0 >0$ throughout this section including (\ref{eqn:evidchisqeq}).  A more compact formulation appears possible.}

\subsection{Choosing the equivalence boundary value $\lambda _0$}\label{sec:choosingbdy}

This section is somewhat more mathematical than the others and so Appendix \ref{app:gloss} of notation and definitions is provided.
 Recall from Section~\ref{sec:eqtestchisq} that  for $\bp ^{(n)} $ in a neighborhood of  $\bu _r$  and  large $n$, the Karl Pearson statistic has an approximate $\chi ^2_{r-1,\lambda _n}$  distribution, where  $\lambda _n =rn d_n^2$,  and
 $d_n=d(\bp ^{(n)} ,\bu_r)$.  It is therefore convenient  to define equivalence to uniformity as $0\leq \lambda < \lambda _0$ for some $\lambda _0$, which amounts to choosing  $0\leq  d < d _0$ for some $d_0$. The choice of $d _0$ generally depends on context, and the practitioner can make such an evaluation in each application. Nevertheless,  a  specific proposal is offered, based on what works effectively in the routine examples to follow in the next section.   {Our starting point} is to  define equivalence to uniformity  by placing a bound of $100k\%$ on the relative distance of each component of $\bp$  from $1/r$.

\begin{definition}\label{def1}
For fixed $0<k\leq 1$ define $\bp$ equivalent to $\bu _r$ if   $M(\bp,\bu_r)\leq M_0=k/r$.
\end{definition}

While Definition~1 is  what we desire to use, as noted above it is mathematically more convenient
 to use the Euclidean metric:

\begin{definition}\label{def2}
For fixed $0<k\leq 1$ define $\bp$ equivalent to $\bu _r$ if   $d(\bp,\bu_r)\leq d_0=k/\sqrt{r(r-1)}\,.$
\end{definition}

\begin{figure}[t!]
\centering
\includegraphics[width =12cm,height=10cm]{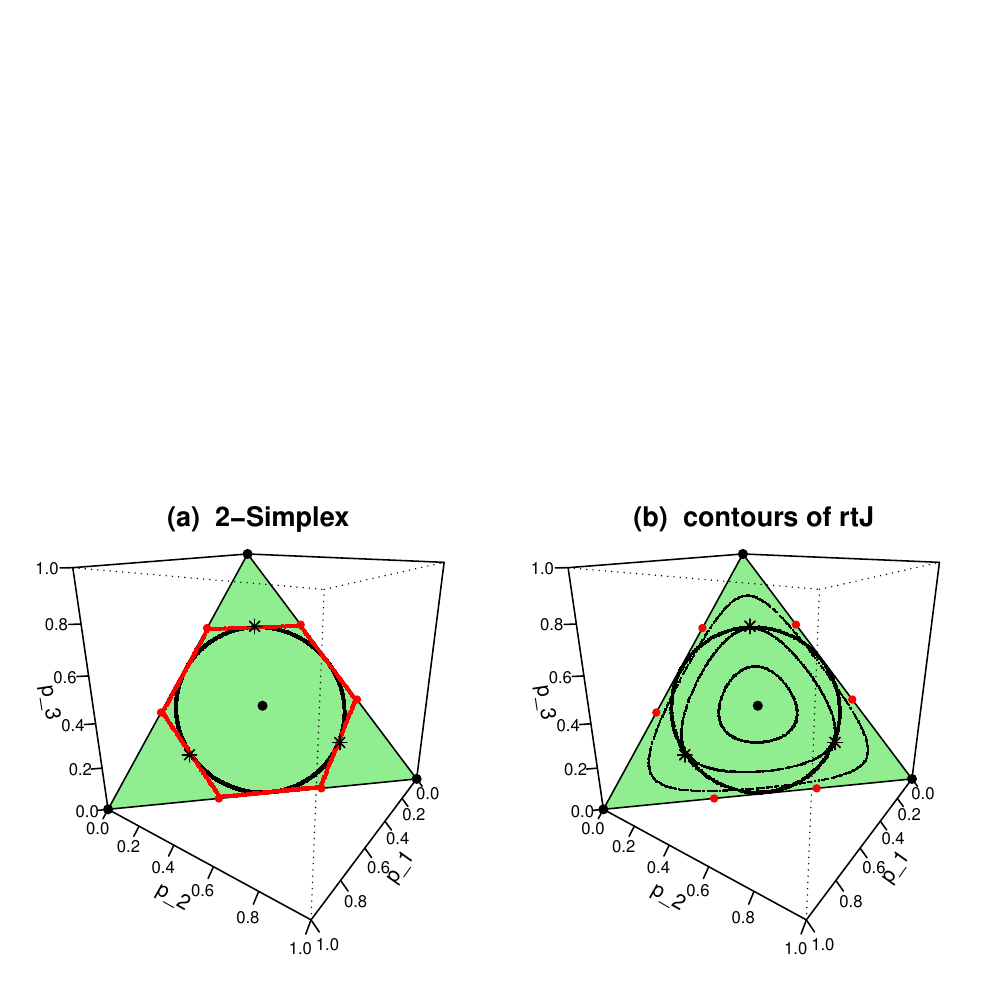}
\caption{ Plot (a) depicts in light-green shading the 2-simplex $\calS _{2}$ in $r=3$  dimensions as the convex hull of its vertices $(1,0,0), (0,1,0)$  and $ (0,0,1) .$
Also shown with boundary in red is $\calC_3(1/3)$, the polytope centered at $\bu _3=(1/3,1/3,1/3)$ and  just contained within $S_2$.  The largest \lq sphere\rq \  within $S_2$ is shown
in black, and has center $\bu _3$ and radius $1/\sqrt{6}$.  Plot (b) shows the same inscribed sphere of the simplex.   Also shown are
three  contours of $\sqrt J\,$, where $J=J(\bp ,\bu _3)$ is the Kullback-Leibler symmetrized divergence of $\bp $ from uniformity $\bu_3$; these contours correspond to
$\sqrt J\, =1/3,2/3,1$ as one moves away from the centre.  See text for interpretation. }  \label{fig6}
\end{figure}

   {To see what compatibility, if any, there is between these definitions, we next examine their geometry for $k=1$.
For $r=2,3,\dots $ denote the probability simplex by  $\calS _{r-1}=\{\bp=(p_1,\dots, p_r):  \hbox { all } p_i\geq 0 , \sum _ip_i=1\}$.
 An example for $r=3$ is the green equilateral triangle   $\calS _{2}$ shown in the perspective plot Figure~\ref{fig6}(a).  Note the origin is hidden behind the triangle. The uniform distribution $\bu_3=(1/3,1/3,1/3)$ is the central point of the triangle.  The points in $S_2$ which are Euclidean distance $d=1/\sqrt 6\,$ from the centre
 are shown as the black circle; it is the inscribed \lq sphere\rq \  within $S_2.$   Also define the polytope centered at $\bu _r$ by
 $\calC_r(M_0)=\{\bp \in \calS _{r-1}:  M(\bp ,\bu _r)\leq  M_0\} .$  An example for $r=3$ is $\calC_3(1/3)$, the hexagon just contained within the triangle   $\calS _{2}$ shown in Figure~\ref{fig6}(a).   For any $r$, the connection between $\calC_r(M_0), $ where $M_0=k/r$ for some $0<k\leq 1$ and its inscribed ball is found in Proposition~\ref{prop1}(a).
  The inscribed ball has  the same centre $\bu _r$ as the polytope and is tangent to each of its faces.  A warning, however;  the similarity of $\calC_r(1/r)$ with its inscribed ball
  observed  for $r=3$ in Figure~\ref{fig6}  rapidly  disappears with increasing $r$. }

  { In   Figure~\ref{fig6}(b)  some contours of $\sqrt J\,$ are shown and suggest that the smallest value of $\sqrt J\,$ is achieved on the inscribed circle
  at the starred point $p^*= (2/3,1/6,1/6)$ and its  permutations.  This is the content of  Proposition~\ref{prop1}(b).  Numerical computations show that
  for points in the inscribed circle of radius $1/\sqrt 6\,$, the minimum value is $\sqrt J\,= 2/3$  while the maximum value  is $\sqrt J\,= 2$, achieved at
  (1/2,1/2,0) and its permutations. }

  {For $r=4$, the reader may verify that $\calC_4(1/4)$ is the octohedron just contained within the tetrahedron $\calS _{3}.$
Also,  the point $\bm {p}^*= \bm {u}_r+  M_0(1,-1/(r-1),\dots, -1/(r-1))$ and its permutations are of interest; an example is Model 7 in Figure~\ref{fig3}.
 These points satisfy $M(\bm {p}^*,\bm {u}_r)=M_0$ and  $d_0=d(\bm {p}^*,\bm {u}_r)= M_0\sqrt {r/(r-1)}\,$ and will play a role in
both parts of Proposition~\ref{prop1} to follow.  }

\begin{proposition}\label{prop1}
 {(a)}  {The inscribed ball of  $\calC_r(1/r)$ is the same as the inscribed ball of $\calS _{r-1}$.  It has center $\bu _r$ and radius $1/\sqrt{r(r-1)}\,.$
The points where this ball just touches  $\calS _{r-1}$  are the $r$ permutations of $\bu _r+\frac {1}{r}(-1,1/(r-1),\dots ,1/(r-1))$.
 The points where it just touches   $\calC_r(1/r)$ are all permutations of  $\bu _r \pm \frac {1}{r}(-1,1/(r-1),\dots ,1/(r-1))$.  For $0<k<1$, it follows that the inscribed ball of $\calC_r(k/r) $ is centered at $\bu _r$ and has radius $k/\sqrt{r(r-1)}\,.$}

\medskip
  {(b)} Let $<d_0<  \sqrt{1-1/r}\,$ be fixed.  To  minimize the Kulback-Leibler divergence
$J(\bu ,\bp) =\sum _{i=1}^r (p_i-1/r)\ln(p_i)$  subject to  the constraints
 $\sum _{i=1}^rp_i=1$ and    $ d(\bp  ,\bu_r)=d_0$ it suffices to take  $\bp^*=\bu _r+d_0 \sqrt{1-1/r}\,(1,-1/(r-1),\dots ,-1/(r-1))$
or a permutation  thereof.
\end{proposition}

\noindent {\bf Proof of Proposition 1(a).}     {$\calS _{r-1}$ is the convex hull of its vertices which are
the $r$ unit vectors $(1,0,\dots ,0),\dots, (0,\dots ,0,1)$  and its centroid is $\bu _r.$ The shortest distance from $\bu _r$ to a point $\bp $ in the boundary of $\calS _{r-1}$  is the distance to one of its facets, say the one opposite vertex $ (1,0,\dots ,0)$.  This facet is the convex hull
of the $r-1$ unit vectors $(0,1,0,\dots ,0),\dots, (0,\dots ,0,1),$  and it has centroid $(0,1/(r-1), \dots ,1/(r-1))$; the distance between the centroids of the simplex
and this facet is then easily computed to be $1/\sqrt{r(r-1)}\,.$  This is the in-radius of the simplex  $\calS _{r-1},$ the radius of its inscribed ball, which touches
the simplex at  $(0,1/(r-1), \dots ,1/(r-1))=\bu _r+\frac {1}{r}(-1,1/(r-1),\dots ,1/(r-1))$ and its  permutations.
The polytope   $\calC _{r}(1/r)$ is the convex hull of  the intersections of hyperplanes that are orthogonal to the line segments joining $\bu _r$ and permutations of $\bu _r  \pm \frac {1}{r}(-1,1/(r-1),\dots ,1/(r-1))$. These points are all equidistant from $\bu _r$ and must lie on  the inscribed ball of $\calC _{r}(1/r)$.  They also include all the points where the simplex meets its inscribed ball.  }

\medskip
\noindent {\bf Proof of Proposition 1(b).}  Use Lagrange multipliers; for details, see Appendix~\ref{app:optim}.

\begin{remark}\label{rem3}
Proposition~\ref{prop1} (a) shows that equivalence of Definition~\ref{def2} is more stringent than Definition~\ref{def1} for the same $k$.
 For the chi-squared statistic, the $d_0$  of Definition~\ref{def2} leads to $\lambda_0 =nr d_0^2=nk^2/(r-1)\,.$   It is shown later that  the choice $k=1/2$ generally provides moderate evidence for equivalence when it should while finding negligible or negative evidence for equivalence under nearby
models; see the fitting of normal and Poisson models  {to data} in Sections~\ref{sec:normality} and \ref{sec:poisson}.
\end{remark}\

\bigskip
\begin{remark}\label{rem4}
Proposition~\ref{prop1} (b) identifies the points $\bp ^*$ on the sphere centered at $\bu _r$ with radius $d_0$ that have the
least  divergence $J(\bp ^*,\bu _r)$, and hence the least expected evidence $\sqrt {J}$ to be found in a
statistic with $\chi ^2_{\nu ,\lambda }$ distribution when the null $\lambda =\lambda _0=nrd_0^2$  and the alternative is $\lambda =0$.
It follows that $\bp ^*$ will be hardest to identify  by a statistical test of these hypotheses.  For from (\ref{evidlevelpower}) and any fixed $\alpha $, the power against an alternative  is monotone  increasing in the expected evidence, which is essentially  $\sqrt {J}$.
\end{remark}

\subsection{Choosing the sample size $n$}\label{sec:choosingn}

This asymptotic mean evidence for equivalence (\ref{eqn:keychisq}) is 0 at the boundary  $\lambda =\lambda _0$ and grows as $\lambda $
   decreases to $0$, where it has a maximum
 {\begin{equation}\label{eqn:maxexpevid}
m_0 = \sqrt {\lambda _0 +\nu /2}-\sqrt {\nu /2}~.
\end{equation}}
   Usually one would want this maximum expected evidence $m_0$ to be at least 3.3, because $T_{\lambda _0,\text {ba}}$ is normal with a standard
   error of  one for estimating its expected value.  For the parameters of Figure~\ref{fig5}, one sees that even if
   $\lambda =0$ (exact uniformity), this expected evidence is  $m_0= \sqrt {12 +5 /2}-\sqrt {5/2}\,=2.2$,  which is between weak and moderate.

In general, in order to obtain a desired maximum expected evidence $m_0$,
solve for $\lambda _0$ in $m_0+\sqrt {\nu/2} =\sqrt {\lambda _0 +\nu /2}$; the result is $\lambda _0 =(m_0+\sqrt {\nu/2})^2-\nu /2.$  The desired value of $n$ is then obtained using the relation $\lambda = nrd^2$ from Section~\ref{sec:eqtestchisq}.
Finally, to ensure the expected number of counts in each of $r$ cells is at least 5,  one needs $n_0\geq 5r$.
Hence the  minimal sample size $n_0$ is

\begin{equation}\label{eqn:samplesize}
n_0=\left \lceil \max \left\{\frac{(m_0+\sqrt {\nu/2})^2-\nu /2}{rd_0^2},\; 5r\right \}\right \rceil  ~,
\end{equation}
where $\lceil x\rceil $ is the smallest integer greater than or equal to $x$.
For equivalence in terms of Definition~\ref{def2} this means that  $d_0=k/\sqrt{r(r-1)}\,;$
examples are in Table~\ref{table2}.

\begin{table}[h!]
  \tbl{Values of sample size $n_0$ defined by  (\ref{eqn:samplesize}) required to achieve maximum expected evidence $m_0$  for
equivalence to uniformity.  The sample sizes are for  the maximum relative error $k=1$;  for smaller $k$ multiply these entries by $1/k^2$.}
 { \begin{tabular}{ccccccccccccc}
  \toprule
                               & & $\qquad m_0\;\backslash \; r$            &   &  &       2         &     3  &  4     &   5       &    6   &  10    &    25   & 100                         \\
     \cmidrule{2-13}
                                && 1.645              &   &  &           10   & 15  &   21     &   30    & 40    &   88   &  339   & 2560       \\
                                & & 3.3                  &   &  &          16   & 35   &  57     &  81    & 107   & 225  &   811   &  5676       \\
                                & & 5                     &   &  &          33   & 70   &  112   &  157  & 205   &  416 &  432    & 9441   \\
 \bottomrule
\end{tabular}}  \label{table2}
\end{table}

In the die example  where $r=6$, to obtain moderate expected evidence $m_0=3.3$ for a die that is perfectly uniform, and equivalence defined by allowing a 50\% relative absolute discrepancy from 1/6, that is with  $k=1/2$ in Definition~\ref{def2},  the sample size must be $n_0= 4\cdot 107=428.$
The reader could introduce another criterion  which is more appropriate for their
application. In particular, for small $r=2$ one would likely demand $k=1/10$   while for $r=100$ one would choose a much larger $k$.

\section{More examples of evidence for goodness of fit}\label{sec:exsevid}
   Next the evidence for normality is found in Section~\ref{sec:normality} and
the evidence for a Poisson model in Section~\ref{sec:poisson}; while important in their own right, they also serve as templates for numerous
other parametric models.  Evidence for uniformity of digits produced by a random number generator and in the decimal digits of $\pi $ are in the online supplementary material.

\subsection{Example 2:   Evidence for normality}\label{sec:normality}

Given a sample of $n$ observations, a  standard approach to chi-squared testing for  normality $N(\mu ,\sigma ^2),$ where both parameters are unknown,  is to first  find  the maximum likelihood estimates $(\bar x ,s_x^2)$ of $(\mu ,\sigma ^2)$ using all the data.  Second, specify $r$ intervals $[\bar x +s_x \Phi ^{-1}((j-1)/r), \bar x +s_x \Phi ^{-1}(j/r)]$, for $j=1,\dots ,r$. If  $\bar x, s_x$ are close to their estimands, these intervals (cells) will have approximately equal probabilities under the model  $N(\mu ,\sigma ^2).$ Third, based on the numbers $\bnu =(\nu_1,\dots ,\nu _r)$ of observations falling in the $r$ intervals,
carry out a test for uniformity or find evidence for uniformity as in earlier sections, with uniformity indicating the normal model is compatible with the data.

 When using this procedure, it is sometimes  recommended   to reduce the \df \;  in the chi-squared approximation to the Karl Pearson statistic by the number of estimated parameters, so in this case of normality $\nu = r-3.$   In fact this modification is quite poor for small $r$, leading to exaggerated significance of tests, as explained in detail by \cite{chern-1954,watson-1957}.  In particular the last author recommends that $r\geq 10 $ if one wants a level-0.05 test to have size between 0.05 and 0.06.    There has also been extensive research on the choice of $r$ to maximize power of the chi-squared test, see the content and references in  \cite[Sections 1.6, 2.14]{G-N-1996} and \cite{quine-1985}. The general consensus is that to have non-trivial power against alternatives $r$ should grow to infinity with $n$ but at a smaller rate, such as $\ln (n)$.   In view of these results, it is suggested here to take
\begin{equation}\label{eqn:rnormal}   r=  \max \{10,\lceil \ln (n)\rceil \}~.
\end{equation}

Table~\ref{table3} summarizes the performance of the above described procedure for finding evidence for normality by listing sample means and standard
deviations $\overline T(s_T)$ of the evidence for uniformity (\ref{eqn:biasadjevideq}) of the selected equiprobable $r$ cells, using Definition~\ref{def2} with  $k=1/2$.

\begin{table}[t!]
\tbl{ For each sample size $n$,  10,000 samples  were generated from each of the normal, logistic and Student-$t$ with 5 \df\; families.
 For each sample the evidence for normality was computed using the method of Section~\ref{sec:normality} with $r$ given by (\ref{eqn:rnormal}) and
 the mean and standard deviation $\overline T(s_T)$ were recorded and listed  in the last three columns below.   The maximum expected evidence when the data are normal is $m_0=\sqrt {\lambda _0 +\nu /2}\, -\sqrt{\nu /2}\;$, where $\nu =r-3$ and $\lambda _0=n /\{4(r-1)\}$.}
 { \begin{tabular}{rrrrrrrrrr}
\toprule
   &  \multicolumn{3}{c}{parameters}               &&       \multicolumn{5}{c}{family}          \\
   \cmidrule{2-4}                            \cmidrule{6-10}
     $n$    & $r\; $        &$\lambda _0\quad $  &  $m_0\;$   &&   $\Phi \qquad$        &&  Logistic\;                && $t_5\qquad$   \\
     100     &      10    &   2.78                &  0.63    &&  0.49(0.90)    &&    0.23(0.92)        &&  $-1.42$(1.05)      \\
     400     &      10    &  11.11              & 1.95      &&  1.90(0.91)    &&    1.04(0.94)        && $ -0.04$(1.16)      \\
   1600     &      10    &  44.44              & 5.05      &&   5.05(0.90)   &&    2.57(0.89)        && $+0.06$(1.22)      \\
   6400     &      10    & 177.78             & 11.59   && 11.61(0.90)   &&    5.46(0.86)        && $+0.18$(1.27)        \\
  25600   &      11     &  640.00            & 23.38   && 23.45(0.90)   &&    9.65(0.86)        &&  $-1.06$(1.22)        \\
$100\cdot 4^5$      &12   &2327.27  &46.17    && 46.23(0.91)   &&  16.79(0.90)        &&  $-4.73$(1.29)    \\
$100\cdot 4^6$ &13 & 8533.33        & 90.17   && 90.24(0.93)   &&  29.05(0.85)       &&   $-13.93$(1.33)    \\
\bottomrule
 \end{tabular}}
   \label{table3}
\end{table}

 The parameters of the approximating chi-square distribution are $\nu =r-3$ and $\lambda_0.$   {The maximum expected evidence  $m_0$ defined by (\ref{eqn:maxexpevid}) }
  is important for it gives the expected evidence if the data are indeed normally distributed, and this $m_0$  is known prior to computing the evidence $T=T_{\lambda _0,ba}(S)$.  When the simulated  data are from any normal distribution, the sample mean of  $T$-values are indeed close to the maximum value $m_0$.

 The logistic distribution is very close to the normal  \cite{staudte2017}, so one would expect the corresponding  entries in the next column to be similar.  However the mean evidences are reduced by more than 50\% for the logistic, so if one obtains strong evidence (say $T\geq 5$) for normality which is much less than $m_0$, one knows that it is really strong evidence for a distribution that is close to the normal.  The same phenomenon occurs for Student-$t$ distributions with \df\; larger than five, while the evidence is negative for \df \,less than 5.  Evidence for normality for the case \df \  equal to 5 is  negligible  or negative as shown in the table.
  This is also the case for  the Laplace and many asymmetric distributions not tabled here.    {For a new look at the closeness of shapes of distributions in a location-scale
  free context, see  \cite{staudte2017,staudte2018}.}

\subsection{Example 3:   Evidence for the Poisson distribution}\label{sec:poisson}
 This example differs from the previous one in that the number of cells $r$ is determined by the data and the procedure for combining tail cells.
 However it turns out that the random $r$ is nearly constant and one can employ the Karl Pearson statistic with $\nu =r-2$ \df .
 A more important difference is that now the probabilities  $\bp^\text{comb} $ on the combined cells are not uniform,  so a new
definition of  equivalence is required.
Recall condition (\ref{eqn:lambda}) of Section~\ref{sec:eqtestchisq} required for  the limiting distribution of the Karl Pearson statistic under alternatives,
where  $\lambda $ is the limiting value of $ \lambda _n=  n\sum _{i=1}^r\frac{(p _i-p_i^{(n)})^2}{p_i}$.  Here the target model after combining cells is $\bp^\text{comb} $, and one can write $\lambda _n $ as $n$ times a weighted average of squared relative discrepancies from $\bp^\text{comb} $ with weights  $\bp^\text{comb}. $

\begin{definition}\label{def3}
Write  $\lambda _n=  n\epsilon ^2$ where $\epsilon ^2=\sum _{i=1}^rp_i^\text{comb}\epsilon _i^2$ and  the relative discrepancies are $\epsilon _i=  \{p _i^\text{comb}-p_i^{(n)}  \}/  p _i^\text{comb} $ for $i=1,\dots ,r$.  To achieve formal consistency with  Definition~\ref{def2} , for fixed $0<k\leq 1$
define equivalence to the Poisson model by $\lambda <\lambda _0$, where $\lambda _0=n\epsilon _0^2$ and $\epsilon _0^2=k^2/(r-1).$
\end{definition}

Given observed counts $\bnu =(\nu_0,\nu _1, \nu _2\dots )$  of integers $0,1,2\dots $whose  size is $\sum _j\nu _j =n$, a standard procedure  in
testing for a Poisson($\mu $) distribution,  $\mu >0$, is firstly, to  find the maximum likelihood estimator $\hat \mu=\sum _jj\,\nu _j /n $ of $\mu $ based on all $n$ observations; and secondly, to combine integers with small probabilities under the
Poisson($\hat \mu$) distribution, so as to leave a consecutive set of $r$ integers (cells) on which to calculate the Karl Pearson statistic.
 These $r$ integers, labelled  $(r_0+1,r_0+2,\dots r_0+r)$,  are chosen so that the expected counts in cells $r_0+1$ and $r_0+r$
each exceed 5 under the Poisson($\hat \mu$) model. The following instructions are followed, with $\mu =\hat \mu$.

\smallskip
 \noindent {\bf Procedure for combining cells:\ } {\em  Let $n$ be fixed and $X\sim $Poisson($\mu$). First define $r_0$ to be the least $k$
such that  $nP_\mu (X\leq k) \geq 5$ and define $p^\text{comb}_{1}=P_\mu (X\leq r_0+1).$
Similarly define $r_0+r$ as the greatest $k$ such that $nP_\mu (X\geq k) \geq 5,$ and let  $p^\text{comb}_{r}=P_\mu (X\geq r_0+r).$
For  the remaining $r-2$ cells let $p^\text{comb}_j=P_\mu (X=r_0+j),  $ for $j=2,\dots, r-1.$  }

Having obtained  $r_0,r$ and the $r$-vector of probabilities $\hat \bp^\text{comb}$  one also needs to find the combined cell counts.  Let  $\bnu ^\text{comb}=(\nu  _1^\text{comb},\dots,\nu  _r^\text{comb})$, where
$\nu  _1^\text{comb}=\sum _{j\leq r_0+1} \nu _j$, $\nu  _r^\text{comb}=\sum _ {j\geq r_0+r} \nu _j$ and $\nu  _j^\text{comb}=\nu _{r_0+j}$ for the remaining $r-2$ combined cells. The Karl Pearson statistic (\ref{eqn:KP}) can then be calculated for the $r$-vectors $\bnu ^\text{comb}$ and
$\hat \bp^\text{comb}$ and the evidence for equivalence from (\ref{eqn:biasadjevideq}).

Definition~\ref{def3}  gives useful results when $k=0.5,$ and $\lambda _0= n/4(r-1)$,  see Table~\ref{table4}.
When $\mu =1$, a sample size of  $n=100$ suffices to achieve weak maximum expected evidence $m_0$
for the Poisson model when it is indeed  Poisson(1),  but as $\mu $ increases one needs larger
$n$ to achieve the same result:  for $\mu =20$ one needs $n=1600$ observations.

Another commonly assumed model for count data $X$ is the negative binomial with parameters $r,p$, where $0<p<1$ and $r>0 $
and this is written $X\sim \hbox{NB}(r,p)$.   It is known that $\mu =\e [X]=r(1-p)/p$ and $\var [X]=r(1-p)/p^2$. Unlike the Poisson($\mu$)
model which has variance equal to the mean, the negative binomial has a larger variance than the mean (sometimes called {\em over-dispersion}).
An alternative parametrization in terms of $\mu $ and $\alpha =1/r$, the {\em dispersion parameter},  is obtained by taking $p=r/(\mu+r)=1/(1+\alpha\mu)$.
Then $\e [X] =\mu $ and  $\var [X] =r(1-p)/p^2=\mu +\alpha \mu ^2.$
The greater the value of $\alpha $, the greater the dispersion.  It is of interest to see how much over-dispersion can be tolerated in finding
evidence for the Poisson model, so the above simulation study was repeated  for the same choices of  $n$ and $\mu $ as before but now sampling from the negative binomial distribution with $\alpha =0.01.$
   Comparing results in Table~\ref{table4}  shows that the evidence for
the Poisson($\mu$) model is basically unchanged when $\mu =1$ or 5; that is, the over-dispersion ($\alpha =0.01$) is not picked up  for small $\mu $ and these sample sizes.    However, for $\mu =10$  the expected maximum evidence $m_0$ under the
Poisson model is not achieved with the expected evidences for the Poisson model only half what is expected.
For $\mu =20$ negligible or even negative evidence for the Poisson model is obtained.
Experimentation shows that if the dispersion parameter is $\alpha =1/\mu $ so that
the variance is twice the mean, then for all values of $n$ and $\mu $ shown in Table~\ref{table4} the evidence for the Poisson model
will be negative.

\begin{remark}\label{rem5}
It is of interest to know how $r$ chosen by the above procedure  depends on $n$ and $\mu $.
For large $\mu $ one can use the normal approximation $P_\mu (X\leq j)\approx \Phi \{(j-\mu)/\sqrt {\mu }\}$ to solve for
$r_0 \approx \sqrt{\mu }\;\Phi ^{-1} (5/n)+\mu $. For the other tail,  solving $\Phi \{(k-\mu)/\sqrt {\mu }\}=1-5/n$
yields $k=r_0+r\approx \sqrt{\mu }\;\Phi ^{-1} (1-5/n)+\mu $. Hence as $n\to \infty $ for fixed $\mu $, and using
the  formula $\Phi ^{-1} (1-1/n)\sim \sqrt {2\ln (n)}\;,$ see \cite[p.109]{Das-2008} or Appendix~\ref{app:das},
$  r=r(n) \approx 2\sqrt{\mu }\;\Phi ^{-1} (1-5/n) \sim  \sqrt{8\,\mu \,\ln (n/5) }\; .$
Thus the standard method of combining cells with low expected values under the estimated Poisson model will lead to $r(n)=O(\sqrt {\ln (n/5) }\;)$ growing slowly with $n$ for fixed moderate to large $\mu $, while $r$ is also increasing in $\sqrt{\mu }\,.$
\end{remark}

\begin{table}
 \tbl{ For each $n$,  20,000 samples of  size $n$  were generated from each of several Poisson
  and negative binomial distributions and the evidence for the model
 computed for each sample  as described in Section~\ref{sec:poisson}.  The  means(standard deviations) of the replicated values are  listed.}
{\begin{tabular}{rcrrrcrr}
  \toprule
      &  \multicolumn{3}{c}{ Poisson($\mu=1$) } &&       \multicolumn{3}{c}{ Poisson($\mu=5$)   }    \\
      \cmidrule{2-4}      \cmidrule{6-8}         $n$    & $\bar r $    & $\bar m_0\qquad$ & $\overline T \qquad  $ &&
                       $\bar r $  &  $\bar m_0\qquad $ & $\overline T \qquad  $ \\
     100     & 4.0   &2.08(0.14)   &  2.17(0.85)   &&    8.1         & 0.82(0.03)   &  0.78(0.91)  \\
     400     & 5.0   &3.93(0.08)   &  4.08(0.89)   &&   10.7        & 1.74(0.09)   &  1.76(0.93)  \\
   1600     & 5.9   &7.72(0.30)   &  7.87(0.96)   &&   13.0        & 3.89(0.00)   & 3.93(0.95)   \\
   6400     & 6.0   &16.53(0.00) &16.70(0.92)   &&   14.0        & 8.90(0.06)   & 9.00(0.95)   \\
\midrule
      &  \multicolumn{3}{c}{ Poisson($\mu=10$) } &&       \multicolumn{3}{c}{ Poisson($\mu=20$)   }    \\
      \cmidrule{2-4}      \cmidrule{6-8}                           $n$    & $\bar r$    & $\bar m_0\qquad$ & $\overline T \qquad  $ &&
  $\bar r$       &  $\bar m_0\qquad $ & $\overline T \qquad  $ \\
   100     & 11.4     & 0.50(0.03)   &  0.45(0.93)   &&   15.7      & 0.31(0.02)   & 0.25(0.94)  \\
   400     & 14.9     & 1.16(0.04)   &  1.15(0.94)   &&   20.8      & 0.74(0.02)   & 0.72(0.95)  \\
   1600   & 17.9     & 2.80(0.05)   &  2.83(0.96)   &&   25.0      & 1.91(0.02)   & 1.92(0.95)   \\
   6400   & 20.0     & 6.65(0.03)   &  6.73(0.96)   &&   29.2      & 4.70(0.06)   & 4.73(0.96)   \\
\midrule \\
      &  \multicolumn{3}{c}{ NB($\mu=1,\alpha=0.01$) } &&       \multicolumn{3}{c}{  NB($\mu=5,\alpha=0.01$)  }    \\
      \cmidrule{2-4}      \cmidrule{6-8}               $n$    & $\bar r$   & $\bar m_0\qquad$ & $\overline T \qquad  $ &&
  $\bar r$     &  $\bar m_0\qquad $ & $\overline T \qquad  $ \\
     100     & 4.0         &  2.08(0.14)   &  2.15(0.85)    &&    8.1   & 0.82(0.04)   & 0.77(0.93)  \\
     400     & 5.0         &  3.93(0.08)   &  4.06(0.90)    &&   10.7  & 1.74(0.09)   &  1.64(0.97)  \\
   1600     & 5.9         &7.73(0.31)      &  7.84(0.98)   &&   13.0  & 3.89(0.00)   & 3.50(1.04)   \\
   6400     & 6.0         &16.53(0.00)   &16.60(0.93)    &&   14.0  & 8.90(0.06)   & 7.61(1.09)   \\
\midrule
      &  \multicolumn{3}{c}{ NB($\mu=10,\alpha=0.01$)  } &&       \multicolumn{3}{c}{ NB($\mu=20,\alpha=0.01$) }    \\
      \cmidrule{2-4}      \cmidrule{6-8}                                                                       $n$    & $\bar r$   & $\bar m_0\qquad$  & $\overline T \qquad  $ &&
  $\bar r$       &  $\bar m_0\qquad $ & $\overline T \qquad  $ \\
   100     & 11.4           &  0.50(0.03)   &  0.45(0.93)   &&   15.7        & 0.31(0.02)   & 0.03(0.98)  \\
   400     & 14.9           & 1.16(0.04)    &  1.15(0.94)   &&   20.8        & 0.74(0.02)   &$ -0.42$(1.13)  \\
   1600   & 17.9           & 2.80(0.05)    &  1.54(1.12)   &&   25.0        & 1.91(0.02)   &$ -1.48$(1.21)   \\
   6400   & 20.0           & 6.65(0.04)    &  3.13(1.13)   &&   29.2        & 4.70(0.06)   & $-3.86$(1.27)   \\
\bottomrule
 \end{tabular}}
\label{table4}
\end{table}

\subsubsection*{ Alpha Particle Emissions Data}
As a specific example, consider  the  Alpha Emissions Data freely available at

 {\tt http://www.randomservices.org/random/}.   It is described there  by
\begin{quote}\small
{\em Americium (atomic number 95) is a synthetic element that is produced as a byproduct in certain nuclear reactions. It was first produced by Glenn Seaborg and his colleagues at the University of California, Berkeley. The isotope americium-241 now has commercial applications in the ionization chambers of smoke detectors. It decays by emission of alpha particles and has a half-life of 432.2 years. In 1966, the statistician J. Berkson studied alpha particle emissions from a sample of americium-241. The table below is a frequency distribution for the number of emissions in 1207 ten-second intervals, and is adapted from data in Rice.}
 \end{quote}
The table is omitted here and the references are \cite{berk-1966} and \cite{rice1993}.
The question for us is whether the observed counts are consistent with a Poisson($\mu $) model.
These counts are $\bnu  = (1,4,13,28,56,105,126,146,164,161,123,101,74,53,23,15,9,3,1,1)$ on 	$0,1,2,\dots ,19.$	 																																																																																																																																																																																																																																																																																																																																																																																																																																																																						
The total sample size is $n=1207$ and the sample mean and variance of these data are $\bar x= 8.367$ and $s_x^2=8.469$, so there is no reason to suspect under- or over-dispersion.  After combining cells 0-2 and also 17-19 using the standard procedure, there are $r=16$ remaining upon which to calculate
the Karl Pearson statistic which is $8.95$ so the traditional test for lack of fit has p-value $P(\chi ^2_{12}\geq 8.95)$=0.84.
But how much evidence for the Poisson($\mu $) model is there in these data?   Equivalence to this model according to  Definition~\ref{def3}
with $k=1/2$ leads  to $\lambda _0=20.117$
and a maximum expected evidence $m_0=2.56$ for these parameters.   In fact the actual evidence is $  T_{\lambda _0,\text {ba}}(S) =3.53\pm 1$ which is  moderate evidence for the Poisson model.

\section{Summary and further research questions}\label{sec:summary}

 {At the very least, we recommend that traditional tests for lack of fit  be replaced by equivalence tests for goodness of fit,
 even though this means one  needs to  specify what is meant by equivalence to the desired model.  We further recommend that such an equivalence test be  supplemented  by, or even replaced by finding the evidence $T$ for the alternative hypothesis of equivalence.  This evidence for equivalence is informative and easy to interpret, and should be examined
  prior to deciding whether to adopt a model based on a given data set. The  many examples in the text and supplementary materials suggest there is more information in the Karl Pearson statistic than  routine testing with it  will reveal.}

The expected evidence  is found numerically to be very close to the square root of the symmetrized Kullback-Leibler divergence between $\chi ^2_{\nu ,\lambda _0}$ and $\chi ^2_{\nu ,\lambda },$ for $0\leq \lambda \leq \lambda _0$. A proof of this approximation with error term, as was found for exponential families in \cite{M-S-2012}, is not yet available. In fact, a more general approximation theorem of this type involving a much wider class of distributions is probably true.

The maximum expected evidence $m_0$ occurs if the desired model  is true, is easily computed and guides one in assessing an observed $T$.
This methodology was illustrated by finding evidence in the Karl Pearson statistic for  equivalence to the normal model when it is actually a good one  and also when the data better fits a nearby one.   Similarly in the case of discrete data where a Poisson model is at question and reduction to an equivalence test to uniformity is not available, one can still compute the evidence for equivalence to the model. Further experimentation could reveal how the
expected evidence depends on the negative binomial parameters $\alpha $, $\mu $ and $n$ when over-dispersion is present.  Evidence for
many other models is easily found by slight modification of these two examples.

A source of difficulty arises in choosing $\lambda _0=nr\,d_0^2$ when one wants to define equivalence to uniformity in terms of the sup metric
$M(\bp,\bu _r)$ instead of $d=d(\bp,\bu _r)$, as discussed in Section~\ref{sec:choosingbdy}.  Instead, one could evaluate the statistic  $M(\hat \bp,\bu _r)$, where $\hat \bp$ is the  maximum likelihood estimator of  $\bp $,   as the basis of evidence for equivalence to uniformity.

Straightforward applications of this methodology to other contexts where the Karl Pearson statistic is routinely employed, such as
tests for independence in contingency tables, are clearly possible.  They would require one to specify what is meant by equivalence to independence.

In an ANOVA comparison of possibly different normal populations the non-central $F$ distribution arises, and a more ambitious project would be to
define equivalence  in terms of its \ncp\, and smoothly extend the classical \vst \, of the $F$ statistic as carried out here for the non-central chi-squared
distribution. A start on this project is made in \cite{M-S-2016}, where a linear extension of  the \vst \, was proposed.

 {
\section*{Acknowledgements}  The author thanks Natalie Karavarsamis for helpful comments on an early version of this manuscript
and to  Yuri Nikolayevsky,   Robin Hill and  Luke Prendergast  for discussions and suggestions on a later draft.  In addition, the author is indebted to
the referee and Editors for their challenging and probing questions which improved the content and presentation of the text.}

\bibliographystyle{plain}
\bibliography{staudtebib}

\begin{thebibliography}{10}
\bibitem{KP-1900}
K.~Pearson.
\newblock On the criterion that a given system of deviations from the probable
  in the case of a correlated system of variables is such that it can be
  reasonably supposed to have arisen from random sampling.
\newblock {\em {Philisophical Magazine}}, 50:157, 1900.
\end{thebibliography}
\vspace{5cm}

\section{Appendix}\label{sec:appendix}

\subsection{Brief summary, glossary and  definitions for Section~\ref{sec:choosingbdy}}\label{app:gloss}

The main result from Sections~\ref{sec:introd}--\ref{sec:equivtest}  is that the Karl Pearson statistic $S$, when
having an approximate $\chi ^2_{\nu,\lambda }$ distribution,  can be monotonically transformed via (\ref{eqn:evidchisqeq}),(\ref{eqn:biasadjevideq}) to
 $T_{\lambda _0,\text {ba}}(S) $, called the {\em evidence for equivalence}.
 Its distribution is  approximately normal with variance one and asymptotic mean function
 $\key _{\nu,\lambda _0}(\lambda )\doteq  \e  _{\nu,\lambda }[T_{\lambda _0,\text {ba}}(S)], $ see below.   Here $\lambda \geq 0 $ is unknown, while  $\nu $ is known and depends on $r$.

 The {\em equivalence boundary} $\lambda _0>0$ is chosen in advance to divide the parameter space:
\begin{description}
\item[\qquad \quad \;$\lambda  _0\leq \lambda $\quad \ ]  null hypothesis of non-equivalence
\item[$0 \leq \lambda  <  \lambda _0$\quad \qquad ]  alternative hypothesis of equivalence
\item[ ]{ }
\item[\quad $\key _{\nu,\lambda _0}(\lambda )$\qquad ] $\key _{\nu,\lambda _0}(\lambda )= \sqrt {\lambda _0+\nu /2}-\sqrt {\lambda +\nu /2}$,  {\em asymptotic mean function (\ref{eqn:keychisq})}
 \item[\quad $m _0$\qquad \qquad \qquad ]  $m _0= \sqrt {\lambda _0 +\nu /2}-\sqrt {\nu /2},$ {\em the maximum expected evidence (\ref{eqn:maxexpevid})}
   \end{description}
 The notation in Sections~\ref{sec:introd}--\ref{sec:equivtest} is standard and derives from  the multinomial distribution  with parameters $n$ and $\bp _r=(p_1,\dots ,p_r).$ In Section~\ref{sec:equivevid} each $\bp _r$ is identified with a point in the $(r-1)$-simplex lying in $r$ dimensions. For example, the uniform distribution $\bu _r=(1/r,\dots ,1/r)$   corresponds to the centroid, or average of the vertices, of this simplex.  The geometry of some neighborhoods of $\bu _r $ are then examined to see what might be useful in deciding what is meant by \lq equivalence to uniformity\rq.\
\begin{description}
  \item[$r$\qquad ]  {a fixed number $r$ of cells, where $r\geq 2$ }
  \item[$\bp_r$\ \quad ] {$(p_1,\dots, p_r):  \hbox { all } p_i\geq 0 , \sum _{i=1}^r p_i=1$,  a categorical model with labels $1,2,\dots ,r$}
  \item[$\calS _{r-1}$\quad  ] {$\calS _{r-1}=\{\bp_r |\hbox { all } p_i\geq 0 , \sum _{i=1}^r p_i=1\}$ {\em probability simplex} in $r$ dimensions}
  \end{description}

\clearpage\newpage

The {\em Euclidean distance} between two vectors $\bx=(x_1,\dots ,x_r) ,\by=(y_1,\dots, y_r) $ in $r$ dimensions is  $d(\bx ,\by )=\{ \sum _{i=1}^r(x_i-y_i)^2\}^{1/2}$ and the {\em supremum distance} is  $M(\bx ,\by)=\max _{1\leq i\leq r}\{|x_i-y_i| \}.$
Examples of $d=d(\bp ,\bu _6)$ and $M=M(\bp ,\bu _6)$ for the Models 1--9 of Figure~\ref{fig3} are listed in Table~\ref{table1} of Section~\ref{sec:eqtestwellek}.

\bigskip
It is mathematically convenient to put a bound $d_0$ on $d=d(\bp ,\bu _r)$,  because then the important parameters $\lambda _0= nrd_0^2$ and $m_0$ are determined.  However, for many researchers,   when defining equivalence to uniformity, it is important to place a bound $M_0=k/r$ on the maximum relative discrepancy  $M(\bp ,\bu _r)=\max _{1\leq i\leq r}\{|p_i-1/r| \}$  between each cell probability $p_i$ and the uniform value $1/r$.
\begin{description}
  \item[$d_0$\quad \qquad \qquad ]   An upper  bound on the Euclidean distance of $\bp _r$ from $\bu _r$.
  \item[$M_0=k/r$\qquad ]   An upper  bound on the sup distance of $\bp _r$ from $\bu _r$.
  \item[$k$\qquad \qquad \quad \  ]   {The relative sup distance one is willing to accept in defining  }
 \item[   $\qquad \qquad \qquad $]  {equivalence  to uniformity. Usually $0\leq k\leq 1$.}
 \item [$\lambda = nrd^2$\qquad ]  {Relation between $d=d(\bp _r,\bu _r)$ and $\lambda $ based on Equation~(\ref{eqn:lambda}).  }
     \end{description}
\noindent Proposition~\ref{prop1}(a) shows that taking $d_0=k/\sqrt {r(r-1)}\,$ will guarantee $M\leq M_0=k/r$.
If $d_0=k/\sqrt {r(r-1)}\,$, then $\lambda _0=nk^2/(r-1).$

\bigskip
One can define neighborhoods of the centroid in the simplex $\calS _{r-1}$:
\begin{description}
  \item[$\calB_r(d_0)$\quad   ] { $\{\bp \in \calS _{r-1}|  d(\bp ,\bu_r)\leq d_0  \} $  } the Euclidean ball of radius $d_0$ centered at $\bu _r$
  \item[$\calC_r(M_0)$\quad ] {$\{\bp \in \calS _{r-1}|  M(\bp ,\bu_r)\leq M_0  \} $ } the polytope  centered at $\bu _r$ with distance to each of its facets equal to $M_0$.
  \end{description}
Examples of $\calS _{2}$, $\calB_3(1/\sqrt{6})$ and $\calC_3(1/3)$ are shown in Figure~\ref{fig6}(a).

\bigskip

The {\em Kullback-Leibler information} $I(0:1)$ in $X\sim f_0$ for discrimination between two probability distributions $f_0,f_1$
is defined in \cite{Kullback-1968} as $I(0:1) = \e _0[\ln (f_0(X)/f_1(X))] .$   It is also called the directed divergence.
The {\em symmetrized Kullback-Leibler divergence}, or \kld, is  $J(0,1)=I(0:1)+I(1:0).$
 An example is calculated for two multinomial distributions in Section~\ref{app:kld}. The special case of $n=1$ gives the symmetrized divergence $ J(\bp,\bu_r)$ of any $\bp $ in the $(r-1)$-simplex $\calS_{r-1}$ from its centroid $\bu _r$.\:

\begin{description}
 \item [$ J(\bp,\bu_r)$\qquad  ] {$J(\bp,\bu_r)=\sum _{i=1}^r (p_i-r^{-1}) \ln(p_i)$.}
  \end{description}

Examples of $J=J(\bp ,\bu _6)$ for the Models 1--9 of Figure~\ref{fig3} are listed in Table~\ref{table1} of Section~\ref{sec:eqtestwellek}.
 Some contours of $\sqrt J=\sqrt {J(\bp,\bu_3)}\,$ in $\calS _2$ are shown in Figure~\ref{fig6}(b).
These plots suggest that this semi-metric might lead to natural convex neighborhoods of the centroid $\bu _r$, but we do not
explore this further here.

\clearpage
\newpage

\subsection{Finding  the KLD for multinomial distributions}\label{app:kld}

Let   $\bp=(p_1,p_2,...,p_r)$ and $\bq=(q_1,q_2,...,q_r)$ .   Observe $\bx=(x_1,...,x_r) $ with multinomial probability function
$ f_0(\bx)= c_n \prod _i p_i^{x_i}$ where $c_n=n!/\prod _i x_i!$ and $\sum _ix_i=n.$ Let $f_1$ denote the model with $\bp$ replaced by $\bq.$
First  find $I(0:1)=\e _0 [\ln (f_0(X)/f_1(X))]$.
\begin{eqnarray}\nonumber
I(0,1)   &=& \e _0\left  [\ln \left \{\prod _i(p_i/q_i)^{X_i}\right \} \right ]   \\  \nonumber
            &=&   n  \left  [\sum _i p_i\left \{\ln (p_i/q_i )\right \}   \right ]   ~.
\end{eqnarray}
By symmetry, $I(1,0)= n  \left  [\sum _i q_i\left \{\ln (q_i/p_i )\right \}   \right ]  $.
Hence the sum $J(0,1) =I(0:1)+I(1:0)$ is
\begin{equation}\label{eqn:Jmultinom}
  J(\bp,\bq)=n\left \{ \sum _i  (p_i-q_i) \ln(p_i)   +  (q_i-p_i) \ln(q_i) \right \}=n\left \{ \sum _i  (p_i-q_i) \ln(p_i/q_i)\right \}   ~ .
\end{equation}

A special case is when  $\bp =\bu = (1/r,...,1/r)$  and to obtain
\[ J(\bp,\bq)=n\left \{ \sum _i  (r^{-1}-q_i) \ln(1/(rq_i))\right \} = n\left \{ \sum _i (q_i-r^{-1}) \ln(q_i)\right \}~.\]

\subsection{Bias in  $T$ for non-uniformity}\label{app:biaspos}
Recall from  Equation~\ref{eqn:evidpos} of Section~\ref{sec:evidpos} that $T=h(S)$, where $h(s)$ is composed of two parts, $h(s)=h_0(s)=\sqrt {2s}-\sqrt{2\nu }\,$ for $s<\nu $  and  $h(s)=h_1(s)=\sqrt {s-\nu/2}-\sqrt{\nu /2} \,$ for $s\geq \nu.$ For any twice differentiable $h(s)$ one has the approximation $\e [h(S)]  \doteq h(\e [S])  +h''(\e [S]) \,\var [S]/2 $, see \cite[Eq.17.1]{KMS-2008} for example.  While the $h(s)$ under consdieration here
is continuously differentiable for all $s>0$, its second derivative is discontinuous at $s=\nu$, which complicates a careful analysis of bias.

First consider the lead term in the expansion for $\e [h(S)] $; it is composed of $h_0(\e  _{\nu,\lambda}[S])=\sqrt {2(\lambda +\nu)}-\sqrt{2\nu })\, $ for
$\lambda < \nu ,$ and $h_1(\e  _{\nu,\lambda}[S])= \sqrt {\lambda +\nu/2}-\sqrt{\nu /2}$ for $\lambda \geq \nu $. Thus $h(\e [S])$ is discontinuous at $\lambda =\nu .$    Further,  $h_1(\e  _{\nu,\lambda}[S]) $  can also be defined for all $0\leq \lambda < \nu $,  and  over this domain   the difference $h_0(\e  _{\nu,\lambda}[S]) -h_1(\e  _{\nu,\lambda}[S])$ is small, in fact less than $\sqrt \nu \,/17.$  Therefore, first  replace $h_0(\e  _{\nu,\lambda}[S])$ by this extended $h_1$ to obtain a smooth asymptotic mean, displayed in (\ref{eqn:ETpos}) .
 This choice is ultimately justified  by simulation studies, which show that  $\e  _{\nu,\lambda} [T]$ is indeed close to  (\ref{eqn:ETpos}) for $\nu \geq 1$
and $\lambda \geq 0$, and can be made closer by a simple bias adjustment. To see what this bias adjustment might be,  we proceed formally
finding:
  $h_0'(s)=(2s)^{-1/2} $ and  $h_0''(s)=-(2s)^{-3/2} $;
and  $h_1'(s)=(s-\nu /2)^{-1/2} /2$ and  $h_1''(s)=-(s-\nu/2)^{-3/2}/4$.
Therefore
\begin{eqnarray*}
\bias _{\nu,\lambda} [T(S)] &=&\left\{  \begin{array}{ll}\e_{\nu,\lambda}[T(S)]-(\sqrt {2(\lambda +\nu)}-\sqrt{2\nu })\,\\ \e_{\nu,\lambda}[T(S)]-(\sqrt {\lambda +\nu/2}-\sqrt{\nu /2})\end{array} \right. \\
&=& \left\{  \begin{array}{ll}  - \frac {2\nu +4\lambda }{2\{2(\nu+\lambda )\}^{3/2} }    , & \hbox{ for $0\leq \lambda < \nu $ ;} \\     -\frac{1}{8} \frac {(2\nu +4\lambda )}{\{\lambda +\nu/2 \}^{3/2} }   , & \hbox{ for $\nu \leq \lambda $ .} \end{array} \right.
\end{eqnarray*}

\[\bias _{\nu,\nu} [T(S)] =\left\{  \begin{array}{ll}  - \frac {6\nu }{2\{4\nu \}^{3/2} } = -  \frac {3 }{8\sqrt \nu}\,   -\frac{0.375}{\sqrt {\nu}\, }  , & \hbox{ for $ \lambda \uparrow \nu $ ;} \\
           -\frac{1}{8} \frac {6\nu}{\{(3\nu)/2 \}^{3/2} }=  -\frac{1}{\sqrt {6\nu}\, } \approx   -\frac{0.408}{\sqrt {\nu}\, }  , & \hbox{ for $ \lambda \downarrow \nu$ .}
\end{array} \right. ~ .\]
By means of simulations it was  found that adding a term such as $\frac{0.408}{\sqrt {\nu}\, }$ to  $T(S)$ did reduce its   bias
at $\lambda =\nu$, but overcompensated at $\lambda =0$ .  After further experimentation, it
 was decided that $  \frac{0.2}{\sqrt {\nu}\, }$ was a useful compromise and so  {version (\ref{eqn:biasadjevidpos}) was adopted.}

\subsection{Bias in  $T_{\lambda _0}$ for equivalence}\label{app:bias0}
Recall from  (\ref{eqn:evidchisqeq}) that  $T_{\lambda _0}= h(S)$ where $h$ is composed of two parts.
Proceeding as in Section~\ref{app:biaspos} it was found that the lead term in $\e [h(S)]  \doteq h(\e [S])  +h''(\e [S]) \,\var [S]/2 $ is also composed of two parts
and can be replaced by extending the second part to the domain $0\leq \lambda <\nu $; this is Equation~(\ref{eqn:keychisq}).
The bias   in $T_{\lambda _0}(S)$ for   (\ref{eqn:keychisq})  is obtained from
For $h_+(s)=-\sqrt {s-\nu/2}$,  $h_+'(s)=-(s-\nu /2)^{-1/2} /2$ and  $h_+''(s)= (s-\nu/2)^{-3/2}/4 $, so the bias term is
\[\bias _+ = \frac{1}{8} \frac {(2\nu +4\lambda )}{\{\lambda +\nu/2 \}^{3/2} } =  \frac {1}{2\{\lambda +\nu/2 \}^{1/2} }       ~ .\]
While this clearly depends on the unknown $\lambda $ it suffices to correct   $T_{\lambda _0}(S)$   at the point $\lambda =\lambda _0$
to obtain a nearly unbiased estimator for the parameters of interest.

 {These bias adjustments Sections~\ref{app:biaspos}, \ref{app:bias0} assume $\nu $ is fixed, but as an anonymous referee has pointed out,  it is of interest to
know how $\nu $ can grow with $n$ in these results.}

\subsection{DasGupta's formula}\label{app:das}

DasGupta \cite[p.109]{Das-2008} states that as $n\to \infty $, $\Phi ^{-1}(1-1/n)\sim \sqrt{2\,\ln(n)}\,$ and that this follows from
the well known asymptotic result $1-\Phi (t) \sim \varphi(t)/t $ as $t \to \infty $.
To verify this claim, let $r_n=\Phi ^{-1}(1-1/n)/ \sqrt{2\,\ln(n)}\,$; then,
 since all $r_n>0$ it suffices to show $r_n^2\to 1$ as $n\to \infty $.  Substituting     $t=  \Phi ^{-1}(1-1/n) $ in $r_n^2$ and applying  L'Hospital 's Rule
yields:
\[ r_n^2\;=\; \frac  { -t^2 }{2\,\ln(1-\Phi(t))}\; \sim \;\frac{t(1-\Phi(t)) }{\varphi (t)}\; \to \;1 ~. \]

\subsection{Proof of Proposition 1(b)}\label{app:optim}

 $J=J(\bu ,\bp) =\sum _{i=1}^r (p_i-1/r)\ln(p_i)$  and $d^2=d^2(\bu ,\bp)=\sum _{i=1}^rp_i^2 -1/r. $
 Suppose $d>0$.  What $\bp $ minimizes $J$ subject to  $\sum _{i=1}^rp_i=1$ and    $\sum _{i=1}^rp_i^2=d^2+1/r$ ?
The Lagrangian for maximizing $-J$ is:
\[L(\bp) =\sum _{i=1}^r (1/r-p_i)\ln(p_i) -\lambda \left(\sum _{i=1}^rp_i-1\right ) -\mu \left (\sum _{i=1}^rp_i^2-d^2-1/r\right ) ~.\]
Setting its partial deratives with respect to unknown variables to 0 yields:

\begin{eqnarray*}
0= \frac{\partial}{\partial p_i} L(\bp)&=& \frac{1}{rp_i}-1-\ln(p_i)-\lambda -2\,\mu \,p_i  \text{\qquad for } i=1,2,\dots,r ~\\
0= \frac{\partial}{\partial \lambda } L(\bp)&= & 1-\sum _{i=1}^rp_i    \\
0= \frac{\partial}{\partial \mu } L(\bp)&=&  d^2+1/r- \sum _{i=1}^rp_i^2  \\
\end{eqnarray*}
The first $r$ equations suggest that possibly all $p_i=1/r$ but that violates the last equation because $d>0$.
So suppose $p_2=p_3=\dots =p_r=p$ for some $0<p<1/(r-1)$ and then the second last equation gives $p_1=1-(r-1)p$. The last equation
yields $p$ as a function of $d$, after solving $(r-1) p^2+    1-2(r-1)p +(r-1)^2p^2=d^2+1/r  $ or $r(r-1) p^2 -2(r-1)p  +1=1/r+d^2$
or $p^2 -2p/r  +c= 0 $, where $c=  (1-1/r-d^2 )/(r(r-1)) =1/r^2-d^2/(r(r-1)).$
Hence
\[ p =1/r\pm \sqrt{1/r^2-c}\,=1/r\pm d/\sqrt{r(r-1)} \]
For the minus sign choice, the requirement $0<p=1/r-d/\sqrt{r(r-1)}$ means $d< \sqrt{(r-1)/r}$ for there to be a solution.
This leads to  $p_1=1-(r-1)p=1/r+d\sqrt{1-1/r}\,.    $
One could solve  for the two unknowns $\lambda ,\mu $ in the equations:
\begin{eqnarray*}
 0&=&\frac{1}{rp}-1-\ln(p)-\lambda -2\,\mu \,p   \\
 0&=&\frac{1}{rp_1}-1-\ln(p_1)-\lambda -2\,\mu \,p_1   \\
 \end{eqnarray*}
Therefore a  possible solution for minimizing $J=J(\bu ,\bp)$ subject to the constraints  is given by
 $\bp^*=\bu _r+d \sqrt{1-1/r}\,(1,-1/(r-1),\dots ,-1/(r-1)).$ One can then check with examples that the solution
is indeed a minimizer.
\end{document}